%% file: pinning.tex
\documentclass[11pt,a4paper]{amsart}
\usepackage[T1]{fontenc}
\usepackage{amssymb}
\usepackage{enumerate}
\usepackage{ifthen}
\usepackage{tabularx}
\usepackage{graphicx}
\usepackage{hypbmsec}
\usepackage{calc}
\usepackage[backref]{hyperref}

\makeatletter
 \theoremstyle{plain}    
 \newtheorem{thm}{Theorem}[section]
 \numberwithin{equation}{section} 
 \numberwithin{figure}{section} 
 \usepackage{verbatim}
 \theoremstyle{plain}
 \theoremstyle{definition}
  \newtheorem{example}[thm]{Example}
 \theoremstyle{remark}    
 \newtheorem{notation}[thm]{Notation} 
 \theoremstyle{remark}    
 
 \theoremstyle{definition}
 \newtheorem{defn}[thm]{Definition}
 \theoremstyle{plain}    
 \newtheorem{prop}[thm]{Proposition} 
 \theoremstyle{plain}    
 \newtheorem{lem}[thm]{Lemma} 
 \theoremstyle{plain}    
 \newtheorem{cor}[thm]{Corollary} 
 \theoremstyle{definition}
 \newtheorem*{defn*}{Definition}
 \theoremstyle{plain}    
 \newtheorem*{prop*}{Proposition} 
 \newtheorem{question}{Question}
 \theoremstyle{remark}
 \newtheorem{remark}[thm]{Remark}
 
\include{macros}

\newcommand{\aug}[2]{{#1}\to{#2}}
\newcommand{\naug}[2]{{#1}\nrightarrow{#2}}
\newcommand{\same}[2]{{#1}\leftrightarrow{#2}}

\newcommand{\parms}{\Upsilon,\Theta}
\newcommand{\edoparms}{\Lambda,\Pi}

\makeatother

\begin{document}

\title{Pinning quasi orders with their endomorphisms}
\author{James Hirschorn}
\address{Graduate School of Science and Technology, Kobe University, Japan}
\email{\href{mailto:j_hirschorn@yahoo.com}{j\_hirschorn@yahoo.com}}
\thanks{The author acknowledges the generous support of the Japanese Society for the
  Promotion of Science (JSPS Fellowship for Foreign Researchers, ID\# P04301).}
\urladdr{\href{http://www.logic.univie.ac.at/~hirschor/}{http://www.logic.univie.ac.at/\textasciitilde hirschor/}}
\date{October 3, 2006}

\begin{abstract}
Some general properties of abstract relations are closely examined. 
These include generalizations of linearity, and properties
based on `pinning' an inequality by a pair of families of endomorphisms.
To each property we try to associate a canonical definition of an augmentation
(or diminishment) that augments (or diminishes) 
any given relation to one satisfying the desired property. 
The motivation behind this study was to identify properties distinguishing between 
the product ordering and the eventual dominance ordering of the irrationals~$\irr$,
and furthermore to identify their relationship 
as a member of a natural class of augmentations.
\end{abstract}

\maketitle
\vspace{-5pt}
\tableofcontents{}
\listoffigures{}

\section{Overview}
\label{sec:overview}

Various general properties of relations, with the emphasis on  quasi orders, are
considered, and for each property a corresponding class of  augmentations (or
diminishments)  is defined which yields a standard  (e.g.~minimum) 
augmentation (or diminishment) 
of an arbitrary relation to one satisfying the specified property. For example,
we generalize the usual notion of linearity as follows. For two families
$\parms\subseteq S^S$ (functions from $S$ into $S$), a relation $(S,\le)$ is
$(\parms)$-linear (definition~\ref{d-1}) if for all $p,q\in S$ either there is a
$\sigma\in\Upsilon$ such that $\sigma(p)\le\sigma(q)$ or there is a $\tau\in\Theta$
such that $\tau(q)\le\tau(p)$. 

Central to this study are properties based on pinning an equality by a pair of
subfamilies of $S^S$. For example, a relation $(S,\le)$ is $(\parms)$-correct 
(definition~\ref{d-3}) if the inequality $\nleq$ is pinned by some member of $\Upsilon$
with respect to the family $\Theta$, by which we mean that for all $p,q\in S$, if
$p\nleq q$ then there is a $\sigma\in\Upsilon$ which pins $\nleq$ with respect to
$\Theta$, i.e.~$\tau\circ\sigma(p)\nleq\tau\circ\sigma(q)$ for all $\tau\in\Theta$. 
In fact, when we restrict our attention to families of endomorphisms, there are
exactly three interesting pinning properties, associated with the inequalities $<$,
$\nleq$ and $\nless$, respectively. 

To each of the new properties introduced here, namely two generalizations of linearity
and three pinning properties, we attempt to associate a canonical definition of an
augmentation (see definition~\ref{d-2}). For example, corresponding to 
$(\parms)$-correctness is the augmentation $\altaltlen{\parms}$ of $(S,\le)$
defined by $p\altaltle q$ if
\begin{equation}
  \label{eq:5}
  \forall\sigma\in\Upsilon\,\exists\tau\in\Theta\spc
  \tau\circ\sigma(p)\le\tau\circ\sigma(q). 
\end{equation}
Imposing some conditions on $\Upsilon$ and $\Theta$ 
we have that $(S,\altaltlen{\parms})$ is a
$(\parms)$-correct augmentation of $(S,\le)$. Moreover, this is the canonical
augmentation for $(\parms)$-correctness. For example, if $\Pi$ is a subsemigroup of
$(\edo(S,\le),\circ)$, the endomorphisms of $(S,\le)$ under composition, then
$\altaltlen{\Pi,\Pi}$ is the minimum augmentation of $\le$ that is $(\Pi,\Pi)$-correct
(corollary~\ref{o-15}). In section~\ref{sec:interrelationships} 
we determine all of the implications between these
properties and their corresponding augmentations; 
these results are summarized in figure~\ref{fig:1}. 

Let $\N=\{0,1,\dots\}$ denote the set of all nonnegative integers. 
We take the \emph{irrationals} to mean the set of all functions from $\N$ into
$\N$---when given the product topology they are 
homeomorphic to the irrationals numbers of the real line (see e.g.~\cite{MR1321597}
or~\cite{irrationals}).
This investigation resulted from a study comparing 
the product order $\le$ on $\irr$ 
of the usual ordering $0<1<\cdots$ of $\N$, i.e.~$x\le y$ if
$x(n)\le y(n)$ for all $n\in\N$, to the \emph{eventual dominance} order $\lefnt$
on $\irr$, i.e.~$x\lefnt y$ if $x(n)\le y(n)$ for all but finitely many $n\in\N$. 
Indeed this paper is the first in a series currently consisting of two papers,
where the next paper~\cite{irrationals} 
is entitled ``Characterizing the quasi ordering of the
irrationals by eventual dominance''. Although much of the discussion will concern 
arbitrary relations, most of our examples are
augmentations of the poset $(\irr,\le)$ and diminishments of the quasi order
$(\irr,\lefnt)$, as they provided the motivation for the abstract
development. The main results obtained here on this comparison are 1) that
$\lefnt$ is the corrective augmentation of $\le$ by the family of all projections
(see definition~\ref{def:proj}) of the
members of $\irr$ onto some infinite set of coordinates (theorem~\ref{u-6}); 
and 2) that $\lefnt$ is the
transitive augmentation (i.e.~transitive closure) of the strictive
augmentation  of $\le$ by the projections (theorem~\ref{x-11}), where
strictness is the property corresponding to pinning the inequality $<$.

\subsubsection*{Conventions}

We should explain our usage of theorem-like assertions. \textbf{Proposition} is used
to indicate a statement which follows directly from the definitions, and either does not
require any proof, or else can be proved in one or two lines; \textbf{Lemma} is used
instead when the proof is any longer; \textbf{Theorem} indicates a result of
distinguished importance, regardless of the length of the proof; and
\textbf{Corollary} is used to indicate a consequence of either a Lemma or a Theorem. 

\section{Quasi order augmentations and diminishments}
\subsection{Terminology}
\label{sec:terminology}

A \emph{quasi order} (also often called a \emph{preorder}) 
is a pair $(O,\le)$ where $\le$ is a reflexive and transitive relation on $O$. 
A \emph{poset} (\emph{partial order}) is a quasi order 
where the relation is also antisymmetric (see section~\ref{sec:antisymmetry}),
and a \emph{strict poset} is a pair $(P,<)$ where $<$ is irreflexive and
transitive. The \emph{complete quasi ordering} of a set $S$ is the quasi order given
by $p\le q$ for all $p,q\in S$. 

\begin{notation}
\label{not:2}
For any quasi order $(O,\le)$ we write $p<q$ for $q$ \emph{strictly
bounds} $p$ in the strict sense, i.e.~$p\le q$ and $q\nleq p$,
or equivalently $[p]\le[q]$ and $[p]\ne[q]$ where $[p]$ denotes
the equivalence class of $p$ in the \emph{antisymmetric quotient}: the poset of 
equivalence classes over the equivalence relation $p\sim q$ if $p\le q$ and $q\le
p$, ordered by $[p]\le [q]$ if $p\le q$. 

More generally, for any symbol of the form $\lenn{\mathrm y}{\mathrm x}$
representing a relation, we let $\lnn{\mathrm y}{\mathrm x}$ denote the relation
satisfying 
\begin{equation}
\label{eq:2}
p \lnn{\mathrm y}{\mathrm x} q\Iff p
\lenn{\mathrm y}{\mathrm x} q\And q\nlenn{\mathrm y}{\mathrm x}p
\end{equation}
for all $p$ and $q$ in the base set. When at least one of $\mathrm x$ or $\mathrm y$
is nonvoid, we let $\eqnn{\mathrm y}{\mathrm x}$ denote the relation satisfying
\begin{equation}
  \label{eq:6}
  p\eqnn{\mathrm y}{\mathrm x} q\Iff p\lenn{\mathrm y}{\mathrm x}q\And
  q\lenn{\mathrm y}{\mathrm x} p
\end{equation}
for all $p,q$ (we  do not want both $\mathrm x$ and $\mathrm y$ void or else
$\eqnn{\mathrm y}{\mathrm x}$ is just the `$=$' symbol). And $\nlenn{\mathrm
  y}{\mathrm x}$, $\nlnn{\mathrm y}{\mathrm x}$ and $\neqnn{\mathrm y}{\mathrm x}$
are the negations of the respective relations. 
\noindent\textbf{\emph{Important.}} Thus in our notation $(O,<)$
is always strict partial order when $(O,\le)$ is a quasi order. This
disagrees with a common usage where $p<q$ means $p\le q$ and $p\ne q$. 
\end{notation}

\begin{example}
\label{x-1}
We defined the eventual dominance relation $\lefnt$ on $\irr$ 
in section~\ref{sec:overview}. 
In accordance with notation~\ref{not:2}, $x\lfnt y$ means $x\lefnt y$ and $y\nlefnt
x$, and thus
\begin{equation}
  \label{eq:7}
  x\lfnt y\Iff x\lefnt y\And x(n)<y(n)\text{ for infinitely many $n\in\N$};
\end{equation}
also according to notation~\ref{not:2}, $x\eqfnt y$ means $x\lefnt y$ and
$y\lefnt x$, and thus $x\eqfnt y$ iff $x(n)=y(n)$ for all but finitely
many~$n\in\N$; and $x\nlefnt y$,  $x\nlfnt y$ and $x\neqfnt y$ are the negations of
$x\lefnt y$, $x\lfnt y$ and $x\eqfnt y$, respectively (e.g.~$x\nlefnt y$ iff
$x(n)>y(n)$ for infinitely many $n$). 

Take note that $x\lfnt y$ has been used differently in the literature with the
meaning $x(n)<y(n)$ for all but finitely many $n$. 
\end{example}

Recall that the class of all binary relations
can be described as the class of all precategories with at most one
arrow from $p$ to $q$ for every pair of objects $(p,q)$, i.e.~this
defines a binary relation by $p\le q$ iff $p\to q$. And the class
of all quasi orders can be described as the class of all categories
with at most one arrow from $p$ to $q$ for every pair of objects $(p,q)$; while the
class of all posets can be described as the class of all categories with at most one
arrow from $p$ to $q$ for every pair of objects $(p,q)$, and with no invertible
arrows besides the identity arrows. 

Letting $\relc$ denote the category of all (small) relations, the arrows or
\emph{homomorphisms} between any two relations $(S,\le)$ and $(R,\altle)$ consist of
all order preserving maps, i.e.~maps $f:S\to R$ such that $p\le q$ implies
$f(p)\altle f(q)$, for all $p,q\in S$. We write $\hom((S,\le),(R,\altle))$ for the
set of all homomorphisms between $(S,\le)$ and $(R,\altle)$. Note that
$\hom((S,\le),(R,\altle))$ can be viewed as the set of all functors from $(S,\le)$
to $(R,\altle)$ when these two are viewed as precategories themselves. 
Then the usual category $\qo$ of all small quasi orders is a full subcategory of
$\relc$, i.e.~the homomorphisms are again the order preserving maps, 
and the category $\po$ of all posets is a full subcategory of $\qo$. 

Homomorphisms from a relation $(S,\le)$ to itself are called \emph{endomorphisms}, 
and we write $\edo(S,\le)$ for the set of all such endomorphisms,
i.e.~$\edo(S,\le)=\hom((S,\le\nobreak),(S,\le))$. We write $\mono((S,\le),(R,\altle))$ and
$\epi((S,\le),(R,\altle))$ for the set of all \emph{monomorphisms} and the set of all
\emph{epimorphisms} from $(S,\le)$ to $(R,\le)$, respectively,
i.e.~$\mono((S,\le),(R,\altle))$ consists of all order preserving injections from
$S$ into $R$, and $\epi((S,\le),(R,\altle))$ consists of all order preserving
surjections from $S$ onto $R$. We denote $\mono((S,\le),(S,\le))$ and
$\epi((S,\le),(S,\le))$ as $\mono(S,\le)$ and $\epi(S,\le)$, respectively. 

\begin{defn}
\label{d-2}
An \emph{augmentation} of a relation $(S,\le)$ 
is a relation $(S,\altaltle)$, with the same base, for which $\altaltle$ is a
superset of $\le$. These are sometimes called \emph{refinements} 
in the literature (e.g.~\cite{MR1902334}). 
And a \emph{diminishment} of a relation is a relation $(S,\altle)$ with the same
base for which $\altle$ is a subset of $\le$, or equivalently, $(S,\le)$ is an
augmentation of $(S,\altle)$. 
Dimishments are sometimes called \emph{weakenings} in the literature 
(e.g.~\cite{MR1808172}). 
By a \emph{quasi order augmentation} of $(S,\le)$ we mean an augmentation of 
$(S,\le)$ which is moreover a quasi order. \emph{Quasi order diminishments}, 
\emph{partial order augmentations} and so forth are defined analogously.
\end{defn}

\begin{notation}
\label{not:3}
We use the notation $\aug{\le}{\altaltle}$ to state that $(S,\altaltle)$ is an
augmentation of $(S,\le)$, or equivalently, that $(S,\le)$ is a diminishment of
$(S,\altaltle)$. Thus $\same{\le}{\altaltle}$ is the same thing as
${\le}={\altaltle}$. 
\end{notation}
\subsection{Abstract augmentations}

Every augmentation can be described in terms of homomorphisms
in the category $\relc$.

\begin{defn}
For a set $S$ and a relation $(R,\altle)$, a functor $f:S\to (R,\altle)$
defines a relation $\len{f}$ on $S$ by\[
p\len{f}q\If f(p)\altle f(q)\textrm{.}\]
If moreover $(S,\le)$ is a relation and $f\in\hom((S,\le),(R,\altle))$ then we call $\len{f}$
a \emph{homomorphic augmentation} of $\le$, because: 
\end{defn}

\begin{prop}
\label{p-16}
$\len{f}$ is an augmentation of $\le$ iff $f\in\hom((S,\le),(R,\altle))$. 
\end{prop}

\noindent In any case, we always have

\begin{prop}
\label{p-2}
$f\in\hom((S,\len{f}),(R,\altle))$.
\end{prop}

Note that

\begin{prop}
\label{p-18}
If $(R,\altle)$ is a quasi order then so is $(S,\len{f})$.
\end{prop}

\begin{prop}
\label{p-19}
If $f$ is an injection and $(R,\altle)$ is a poset, then $(S,\len{f})$
is a poset.
\end{prop}

While every augmentation $\altaltle$ of a relation $(S,\le)$ can
of course be represented as the homomorphic augmentation $\len{i}$
via the inclusion functor $i:(S,\le)\to(S,\altaltle\nobreak)$, the
point here is that homomorphic augmentations sometimes provide 
a nice representation of an augmentation, 
and they also provide a means of constructing
augmentations with useful properties.

\begin{example}
\emph{Restrictive augmentations.} 
Any partial order $(P,\le)$ can be identified with a subset of its
{}``Dedekind'' completion $(\power(P),\subseteq)$:\[
p\mapsto\ds{p}=\{ q\in P:q\le p\}\textrm{,}\]
i.e.~this map is an \emph{embedding} (i.e.~$p\le q$ iff $\ds{p}\subseteq\ds{q}$)
since $(P,\le)$ is a quasi order, and it is injective since $(P,\le)$ is moreover a poset. 
This view leads to a natural way of augmenting a given poset, and more generally a
given quasi order $(O,\le)$.
For a subset $X\subseteq O$, we can define a functor $f_X:(O,\le)\to
(\power(O),\subseteq)$ by 
\begin{equation}
  \label{eq:42}
  f_X(p)=\ds p\cap X,
\end{equation}
i.e.~$f_X\in\hom((O,\le),(\power(O),\subseteq))$. 
Then
\begin{equation}
  \label{eq:43}
  p \len{f_X} q\Iff \ds p\cap X\subseteq\ds q\cap X,
\end{equation}
and by propositions~\ref{p-16} and~\ref{p-18}, $(O,\len{f_X})$ is a quasi order
augmentation of $(O,\le)$. 
Note that the augmentedness increases as $X$ decreases, 
i.e.~$X\subseteq Y$ implies $\aug{\len{f_Y}}{\len{f_X}}$, and $\same{\len{f_O}}{\le}$.
We call $\len{f_X}$ the \emph{restrictive augmentation} of $(O,\le)$ by $X$. 
\end{example}
\subsection{Projections}
\label{sec:projections}

Projections provide us with the fundamental examples of quasi order
homomorphisms for products of quasi orders. They are to be used as parameters for
augmentations. 

For any family $(S_i,\len i)$ ($i\in I$) of relations where $I$ is some index set,
the Cartesian product is viewed as 
a relation where the comparison is made coordinatewise. 
The subcategories $\qo$ and $\po$ of $\relc$ 
are closed under arbitrary Cartesian products. 

\begin{defn}
\label{def:proj}
Suppose $S_{i}$ ($i\in I$) is a family
of sets where $I$ is some index set. For each $j\in I$,
the map $\pi_{j}:\prod_{i\in I}S_{i}\to S_{j}$ defined by
\begin{equation}
\label{eq:9}
\pi_{j}(x)=x(j)
\end{equation}
is an epimorphism called the \emph{projection onto the $j\Th$ coordinate}.
Generalizing to $h:J\to I$ ($J\subseteq I$), the \emph{projection by $h$} is defined by
\begin{equation}
  \label{eq:10}
  \pi_h(x)=x\circ h,
\end{equation}
and thus $\pi_h:\prod_{i\in I}S_i\to\prod_{j\in J}S_{h(j)}$.

Suppose further that the index set $I$ is well ordered. For each $a\subseteq I$, 
we define 
\begin{equation}
  \label{eq:11}
  \pi_a=\pi_{e_a},
\end{equation}
where $e_a$ is the (unique) enumeration of $a$ by the first $\otp(a)$ elements of
$I$. This is useful when dealing with a power, i.e.~$S_i=S$ for all $i$, in which
case $\pi_a:S^I\to S^{\otp(a)}$. 
\end{defn}

\begin{prop}
\label{p-5}
$\pi_h\in\hom\bigl(\prod_{i\in I}S_i,\prod_{j\in J}S_{h(j)}\bigr)$ for all $h:J\to I$. If
moreover, $h$ is an injection, 
then $\pi_h\in\epi\bigl(\prod_{i\in I}S_i,\prod_{j\in J}S_{h(j)}\bigr)$. 
\end{prop}

Thus the `true' projections are given by injections $h$.

\begin{defn}
\label{d-5}
We write $\proj\bigl(\prod_{i\in I}S_i\bigr)=\{\pi_h:h$ injects $I$ into $I\}$, or
just $\proj$ when the intended product is clear.
\end{defn}

\begin{prop}
\label{p-7}
$\pi_J(x)(j)=x(i)$ where $i$ is the $j\Th$ element of $J$. 
\end{prop}

\begin{prop}
\label{p-6}
Suppose $S_i=S$ for all $i\in I$. If $h:I\to I$ then $\pi_h\in\edo(S^I)$, and if
moreover $h$ is a bijection then $\pi_h\in\aut(S^I)$. Thus when $I$ is well ordered and
$J\subseteq I$ has the same order type as~$I$, $\pi_J$ is an endomorphism. 
\end{prop}

Note that $h\mapsto \pi_h$ is contravariant:

\begin{prop}
\label{p-8}
If $g:J\to I$ and $h:K\to J$ \textup($K\subseteq J\subseteq I$\textup), then 
$\pi_{g\circ h}=\pi_h\circ \pi_g$.
\end{prop}
\begin{proof}
For all $x\in\prod_{i\in I}S_i$,
\begin{equation}
  \label{eq:12}
  \pi_{g\circ h}(x)=x\circ g\circ h=\pi_h(x\circ g)=\pi_h\circ \pi_g(x).\qedhere
\end{equation}
\end{proof}

Some notation specific to the order $\irri I$ will be needed later.

\begin{notation}
\label{not:1}
We denote the \emph{support} of a member $x$ of $\irri I$ by
\begin{equation*}
  \supp(x)=\{i\in I:x(i)\ne0\}.
\end{equation*}
We write $\zero_I$, or just $\zero$, for the element with empty support. More
generally, for each $n\in\N$ we let $\mathbf n$ denote the element satisfying
$\mathbf n(i)=n$ for all $i\in I$. And $\chi_J$ denotes the characteristic
function of $J\subseteq I$, and for each $i\in I$, we write $\chi_i$ for
$\chi_{\{i\}}$. Thus $\supp(\chi_i)=\{i\}$ and $\chi_i(i)=1$.  
\end{notation}
\section{Properties of quasi orders}
\label{sec:prop-quasi-orders}

We discuss various properties of abstract relations, and their associated
augmentations or diminishments. The standard relational properties of antisymmetry
(subsection~\ref{sec:antisymmetry}) 
and transitivity (subsection~\ref{sec:transitivity}) are considered. And we consider 
the separativity property (subsection~\ref{sec:separativity}) 
which is commonly mentioned in the context of set
theoretic forcing with a partial order. Then five new properties are introduced. The
two properties $(\parms)$-linearity (subsection~\ref{sec:upsil-theta-line}) and
strict $(\parms)$-linearity (subsection~\ref{sec:strict-upsilon-theta}) are
generalizations of the usual notion of linearity. The three remaining properties,
$(\parms)$-strictness (subsection~\ref{sec:upsil-theta-strictn}),
$(\parms)$-correctness (subsection~\ref{sec:upsil-theta-corr}) and negative
$(\parms)$-strictness (subsection~\ref{sec:negat-upsil-theta}), correspond to
pinning the inequalities $<$, $\nleq$ and $\nless$, respectively. 

\subsubsection*{Pinning}

Now we introduce the notion of pinning.

\begin{defn}
Let $(S,\rel)$ be a relation. For $\Theta\subseteq S^S$ (i.e.~functions from $S$
into~$S$) and $p,q\in S$, 
we say that a function
$\sigma\in S^S$ \emph{pins} the statement $\ulc p\rel q\urc$ with respect to $\Theta$ if
\begin{equation}
  \label{eq:8}
  \tau\circ\sigma(p)\rel\tau\circ\sigma(q)\espc\text{for all $\tau\in\Theta$}.
\end{equation}
And we say that a family $\Upsilon\subseteq S^S$ \emph{pins} the relation $\rel$
with respect to $\Theta$ if for all $p,q\in S$, $p\rel q$ implies there exists
$\sigma\in\Upsilon$ that pins $\ulc p\rel q\urc$ with respect to $\Theta$.
\end{defn}

A given relation $(S,\le)$ induces the inequality relations on the same base set
$S$, i.e.~the relations $(S,<)$, $(S,\nleq)$ and $(S,\nless)$. We shall consider
pinning for these three inequality relations. Note that pinning is symmetric in the
following sense:

\begin{prop}
\label{p-3}
Let $(S,\le)$ be a relation.
If $\Upsilon$ pins the relation $\le$ with respect to~$\Theta$, then it also pins
the relation $\ge$ with respect to $\Theta$. 
\end{prop}

Thus the relations corresponding to the inequalities $>$, $\ngeq$ and $\ngtr$
are all covered by the above relations. Notice that for the given the relation $(S,\le)$ we
have not mentioned pinning $\le$ itself. This is because we are primarily concerned
with families of functions that are endomorphisms, in which case pinning is
automatic:

\begin{prop}
\label{p-4}
Suppose $\Theta\subseteq\edo(S,\le)$. Then every endomorphism $\sigma$ pins every
instance of $p\le q$ with respect to $\Theta$.
\end{prop}
\begin{proof}
Let $\sigma\in\edo(S,\le)$.
If $p\le q$, then for all $\tau\in\Theta$, $\tau\sigma(p)\le\tau\sigma(q)$ because
$\tau\sigma\in\edo(S,\le)$. 
\end{proof}

It should be mentioned that we deem the primary instance of pinning to be the
case where $\Upsilon=\Theta$. Indeed in the next paper of this series this is the
only case considered.

\begin{notation}
The pair $(S,\le)$ will be used to represent an arbitrary relation. We use $(O,\le)$
to represent an arbitrary quasi order, and $(P,\le)$ for an arbitrary
poset. Henceforth, we shall use the parameter pair $(\Upsilon,\Theta)$ to denote a
pair of subsets of~$S^S$, whereas we use $(\Lambda,\Pi)$ to indicate that the
parameters consist of endomorphisms, i.e.~$\Lambda,\Pi\subseteq\edo(S,\le)$. 
\end{notation}
\subsection{Separativity}
\label{sec:separativity}

\begin{defn}
Two elements $p$ and $q$ of a quasi order $(O,\le)$ are \emph{compatible}, 
written $p\compat q$, if they have a common extension $r\le p,q$.
We write $p\incompat q$ for $p$ \emph{incompatible} with $q$ (i.e.~not compatible). 
\end{defn}

\begin{defn}
A quasi order $(O,\le)$ is \emph{separative} if for all $p\nleq q$ in $O$, 
there exists $r\le p$ in $O$ such that  $r\incompat q$. 
\end{defn}

\begin{prop}
Quasi linear orders 
\textup(s.v.~subsection~\textup{\ref{sec:upsil-theta-line}}\textup)
are either complete or nonseparative.
\end{prop}

\begin{prop}
\label{p-21}
Every quasi order with a minimum element is either complete or nonseparative. 
\end{prop}

\begin{example}
For any set $X$, 
$(\power(X)\setminus\{\emptyset\},\subseteq)$ is separative. 
For supposing $a\nsubseteq b$, say $x\in a\setminus b$, 
then the singleton $\{x\}\subseteq a$, and $\{x\}\incompat b$.

By comparison, for any nonempty index set $I$, $(\irri I\setminus\{\zero_I\},\le)$ is
nonseparative. Consider $2\cdot\chi_i\nleq \chi_i$ (see notation~\ref{not:1}). 
\end{example}

Let $(O,\le)$ be a given quasi order. 
Then let $f:(O,\le)\to(\power(O),\subseteq\nobreak)$
be the functor determined by
\begin{equation*}
  f(p)=\{q\in O:q\compat p\}.
\end{equation*}

\begin{defn}
We write $\seple$ for $\len f$ and call it the \emph{separative augmentation} 
of $(O,\le)$, for reasons explained below. 
\end{defn}

\begin{prop}
\label{p-17}
$\seple$ is quasi order augmentation of $\le$ that is separative.
\end{prop}
\begin{proof}
Clearly $p\le q$ implies $f(p)\subseteq f(q)$, and thus $\len f$ is a quasi order
augmentation of $\le$ by propositions~\ref{p-16} and~\ref{p-18}. 
And it is clearly separative. 
\end{proof}

$\seple$ is the minimum separative augmentation in the following sense.

\begin{lem}
\label{p-22}
$(O,\seple)$ is the minimum augmentation of $(O,\le)$ 
that is separative and preserves incompatibility, 
i.e.~if $\altaltle$ is an augmentation of $\le$, $\altaltle$ is separative 
and for all $p,q\in O$, $(O,\le)\models p\incompat q$ 
iff $(O,\altaltle)\models p\incompat q$,
then $\aug{\seple}{\altaltle}$. 
\end{lem}
\begin{proof}
Suppose $p\naltaltle q$. Then there exists $r\altaltle p$ 
which is incompatible with $q$ with respect to $\altaltle$. Since
$\aug{\le}{\altaltle}$, $(O,\le)\models r\incompat q$. 
And since $\altaltle$ preserves incompatibility, $(O,\le)\models r\compat p$.
Thus $p\nlen f q$.
\end{proof}

\begin{cor}
If $(O,\le)$ is separative, then $\same{\seple}{\le}$. 
\end{cor}

\begin{remark}
\label{remark:1}
Note that the \emph{separative quotient} of a quasi order $(O,\le)$ is the set of
equivalence classes over the relation $p\sim_\sep q$ if $f(p)=f(q)$. 
\end{remark}

\begin{example}
Borrowing notation from analysis $c_{00}=\{x\in\irr:\supp(x)$ is finite$\}$, and
writing $c_{00}^+=\irr\setminus c_{00}$, we consider the lattice $(c_{00}^+,\le)$. 
The separative augmentation of this lattice is given by 
\begin{equation}
  \label{eq:44}
  x\seple y\Iff \supp(x)\subseteqfnt \supp(y).
\end{equation}
Also note that the separative quotient is isomorphic to
$((\pnfin)\setminus \{0\},\subseteqfnt)$, where $\Fin$ denotes the ideal of all 
finite subsets of $\N$, $\pnfin$ is the quotient over this ideal, and $\subseteqfnt$
is inclusion modulo finite. 

By comparison, the separative augmentation of $(\finp,\subseteq)$, where $\finp$ is
the coideal $\pN\setminus\Fin$, is given by
$a\seple b$ iff $a\subseteqfnt b$. 
\end{example}

\subsection{Antisymmetry}
\label{sec:antisymmetry}

Recall that a quasi order $(O,\le)$ is \emph{antisymmetric} iff $(O,\le)$ is a poset
iff $p\le q$ and $q\le p$ imply $p=q$ for all $p,q\in O$. 

Note that if $(O,\le)$ is not antisymmetric then neither are any of its
augmentations. 

\begin{defn}
For a given quasi order $(O,\le)$ we define 
\begin{equation*}
  p\asymle q \If p=q\OR p<q
\end{equation*}
(see notation~\ref{not:2}), and call $\asymle$ the \emph{antisymmetric diminishment}.
\end{defn}

\begin{prop}
$(O,\asymle)$ is a partial order diminishment of $(O,\le)$ \textup(and in particular
it is antisymmetric\textup).
\end{prop}

\begin{prop}
If $\le$ is antisymmetric then $\same{\asymle}{\le}$.
\end{prop}

$\asymle$ is not in general a minimal diminishment that is antisymmetric. For
example:

\begin{prop} 
If $(O,\le)$ is a complete quasi order then the asymmetric
diminishment is the equality partial order, i.e.~$\same{\asymle}{=}$. 
\end{prop}

\noindent However, it is in some sense, which we will not make precise, the minimal
definable antisymmetric diminishment. Obtaining a minimal diminishment generally
requires an application of choice. 

\begin{example}
\label{x-8}
The antisymmetric diminishment of $(\irr,\lefnt)$ satisfies
\begin{equation*}
  x \asymle y\Iff x=y\OR x\lfnt y
\end{equation*}
(see example~\ref{x-1}).
\end{example}
\subsection{Transitivity}
\label{sec:transitivity}

Recall that a relation $(S,\le)$ is \emph{transitive} if for all $p,q,r\in S$, $p\le q$ and
$q\le r$ imply $p\le r$. 

\begin{defn}
The \emph{transitive augmentation} of $(S,\le)$, written $\trnle$, is defined by
$p\trnle q$ if there exists a finite sequence $p_0,\dots,p_{n-1}$ in $S$ which forms
a chain from $p$ to $q$, i.e.~%
\begin{equation*}
p=p_0\le p_1\le\cdots\le p_{n-1}=q.   
\end{equation*}
\end{defn}

\begin{prop}
\label{p-31}
$\trnle$ is an augmentation of $\le$ that is transitive.
\end{prop}

This is usually called the \emph{transitive closure} of a relation in the
literature; however, the terminology transitive augmentation matches our present
context. The following basic fact is well known. 

\begin{prop}
\label{p-32}
$\trnle$ is the minimum augmentation of $\le$ that is transitive, i.e.~if
$\altaltle$ is transitive and $\aug{\le}{\altaltle}$ then $\aug{\trnle}{\altaltle}$. 
\end{prop}
    
Note that endomorphisms are preserved under the transitive augmentation.

\begin{prop}
\label{p-37}
$\edo(S,\le)\subseteq\edo(S,\trnle)$.
\end{prop}

\begin{defn}
A \emph{cycle} in a relation $(S,\le)$ means a (finite) cycle in the associated directed
graph, i.e.~a finite $C\subseteq S$ is a cycle iff it is of the form
$C=\{p_0,\dots,p_{n-1}\}$ where $p_0\le p_1\le\cdots \le p_{n-1}\le p_0$. 
We say that the cycle is \emph{bidirectional} if $p_{n-1}\le p_{n-2}\le\cdots\le p_0$.
\end{defn}

\begin{prop}
\label{p-36}
For all $p,q\in S$, $p\trneq q$ \textup(i.e.~$p\trnle q$ and $q\trnle p$, 
see~\eqref{eq:6}\textup) 
iff there exists a cycle containing both $p$ and $q$. 
\end{prop}

We make the following observation concerning when strictly less than is
preserved under taking the transitive augmentation. 

\begin{prop}
\label{p-33}
The following are equivalent\textup{:
 \begin{enumerate}[(a)]
 \item \textit{$p\trnl q$ iff there exists a chain $p_0\le\cdots\le p_{n-1}$ from $p$
     to $q$ such that $p_{i-1}<p_{i}$ for some $i=1,\dots,n-1$, for all $p,q\in S$.}
 \item \textit{$p<q$ implies $p\trnl q$ for all $p,q\in S$.}
 \item \textit{Every cycle in $(S,\le)$ is bidirectional.}
 \end{enumerate}
}
\end{prop}
\begin{proof}
Use proposition~\ref{p-36}.
\end{proof}

\subsection(Linearity){$(\Upsilon,\Theta)$-linearity}
\label{sec:upsil-theta-line}

Recall that a \emph{quasi linear order} is a quasi order $(O,\le)$ that is \emph{linear}:
$p\le q$ or $q\le p$ for all $p,q\in O$. We generalize the notion of linearity as follows.

\begin{defn}
\label{d-1}
Suppose $\Upsilon,\Theta\subseteq S^{S}$. 
We say that $\le$ is $(\Upsilon,\Theta)$\emph{-linear}
if for all $p,q\in S$, \[
\sigma(p)\le\sigma(q)\textrm{ for some $\sigma\in\Upsilon$}\OR\tau(q)\le\tau(p)
\textrm{ for some $\tau\in\Theta$.}\]
In the special case where $\Upsilon=\Theta$, we just say $\Theta$\emph{-linear.}
\end{defn}

\begin{prop}
\label{p-12}
$\le$ is $\Theta$-linear iff for all $p,q\in S$ there is a $\tau\in\Theta$
for which $\tau(p)$ is comparable to $\tau(q)$.
\end{prop}

\begin{example}
\label{x-18}
The partial order $(\irr,\le)$ is $\proj$-linear (see definition~\ref{d-5}), because
for all $x,y\in\irr$, there is an infinite $a\subseteq\N$ such that either $x(n)\le
y(n)$ for all $n\in a$, or $y(n)\le x(n)$ for all $n\in a$, and then
$\pi_a=\pi_{e_a}\in\proj$ and $\pi_a(x)$ is comparable with $\pi_a(y)$. 
\end{example}

\begin{notation}
For $x\in\prod_{i\in I}X_i$ and $y\in\prod_{j\in J}Y_j$ we let $x\bigext y$ denote
the image of $(x,y)$ under the natural association
between $\bigl(\prod_{i\in I}X_i\bigr)\times\bigl(\prod_{j\in J}Y_i\bigr)$ 
and $\prod_{i\in I\djun J}Z_i$ where $Z_i=X_i$ for $i\in I$ and $Z_j=Y_j$ for $j\in J$.
\end{notation}

\begin{example}
\label{x-20}
We define $\Pi_0\subseteq\edo(\irr,\le)$ as the family of all $\sigma:\irr\to\irr$
of the form $\sigma(x)=\pi_h(x)\bigext \rho(x)$ 
where $h:\even\to\even$ is an injection and $\rho\in\sym(\odd)$, 
i.e.~$\rho$ is a permutation of the odd numbers. Note that
$\Pi_0$ is a submonoid of $(\edo(\irr,\le\nolinebreak),\circ)$. 
And $(\irr,\le)$ is not $\Pi_0$-linear. 
\end{example}

Note that $(\Upsilon,\Theta)$-linearity is symmetrical.

\begin{prop}
\label{p-1}
$\le$ is $(\Upsilon,\Theta)$-linear iff it is $(\Theta,\Upsilon)$-linear.
\end{prop}

The strength of linearity decreases with the parameters $\Upsilon$ and $\Theta$.

\begin{prop}
\label{p-20}
If $\Upsilon\subseteq\Upsilon'$ and $\Theta\subseteq\Theta'$ then
$(\Upsilon,\Theta)$-linearity entails $(\Upsilon',\Theta')$-linearity.
\end{prop}

If $\Upsilon$ or $\Theta$
consists of a single element $\sigma$ we write $\sigma$ instead of $\{\sigma\}$.

\begin{prop}
\label{p-13}
$\le$ is $(\id,\Theta)$-linear iff for all $p,q\in S$,\[
p\nleq q\impls\tau(q)\le\tau(p)\textup{ for some $\tau\in\Theta$.}\]
\end{prop}
 
\begin{prop}
A quasi order is a quasi linear order iff it is $\id$-linear. 
\end{prop}

\begin{defn}
The $(\Upsilon,\Theta)$\emph{-linear augmentation} of $\le$ is defined by 
\begin{equation*}
\begin{split}
p\linlen{\Upsilon,\Theta}q\If &p\le q\OR\\
(&\sigma(p)\nleq\sigma(q)\textrm{ for all $\sigma\in\Upsilon$}
\and\tau(q)\nleq\tau(p)\textrm{ for all $\tau\in\Theta$})\textrm{.}
\end{split}
\end{equation*}
Accordingly, the $\Theta$\emph{-linear augmentation} of $\le$ is
given by ${\linlen{\Theta}}={\linlen{\Theta,\Theta}}$.
\end{defn}

\begin{prop}
\label{pro:lin-simple}
$p\linlen{\id,\Theta}q$ iff $p\le q$ or $\tau(q)\nleq\tau(p)$
for all $\tau\in\Theta$.
\end{prop}

\begin{lem}
\label{l-16}
Assume that there exists $\sigma\in\Upsilon$ such that
$\Upsilon\circ\sigma\subseteq\Upsilon$ and 
$\Theta\circ\sigma\subseteq\Theta$
\textup(e.g.~if $\id\in\Upsilon$\textup). 
Then $\linlen{\Upsilon,\Theta}$ is
a $(\Upsilon,\Theta)$-linear augmentation of~$\le$. 
\end{lem}
\begin{proof}
Take any $p,q\in S$. Then $\sigma(p)\linlen{\Upsilon,\Theta}\sigma(q)$
witnesses $(\Upsilon,\Theta)$-linearity, 
and hence we suppose $\sigma(p)\nlinlen{\Upsilon,\Theta}\sigma(q)$.
This implies that either $\pi\sigma(p)\le\pi\sigma(q)$ for some
$\pi\in\Upsilon$ or $\tau\sigma(q)\le\tau\sigma(p)$ for some $\tau\in\Theta$,
and thus $\pi\sigma(p)\linlen{\Upsilon,\Theta}\pi\sigma(q)$ or
$\tau\sigma(q)\linlen{\Upsilon,\Theta}\tau\sigma(p)$. Since $\pi\sigma\in\Upsilon$
and $\tau\sigma\in\Theta$ this verifies $(\Upsilon,\Theta)$-linearity.
\end{proof}

The definition is reasonable in the following sense.

\begin{prop}
\label{p-23}
If $\le$ is $(\Upsilon,\Theta)$-linear then it is equal to its $(\Upsilon,\Theta)$-linear
augmentation.
\end{prop}

And the augmentedness respects the parameters in accordance with
proposition~\ref{p-20}. 

\begin{prop}
\label{p-24}
If $\Upsilon\subseteq\Upsilon'$ and $\Theta\subseteq\Theta'$ then
$\aug{\linlen{\Upsilon',\Theta'}}{\linlen{\Upsilon,\Theta}}$. 
\end{prop}

We have not established here that this is the best possible definition of an 
$(\Upsilon,\Theta)$-linear augmentation. We cannot achieve minimality with one
definition since minimality in general requires the axiom of choice. 
Also note that despite the fact that $(\Upsilon,\Theta)$-linearity is symmetrical in
its two parameters (proposition~\ref{p-1}), the augmentation
$\linlen{\Upsilon,\Theta}$ is not. In examples~\ref{x-9} and~\ref{x-10}, one of the two
possibilities $\linlen{\Upsilon,\Theta}$ yields a natural augmentation of $(\irr,\le)$, 
while the other $\linlen{\Theta,\Upsilon}$ is unnatural and not even transitive. 
A positive
answer to the following question would at least prove that $\linlen{\Theta}$ is the
right definition of a $\Theta$-linear augmentation in the case $\Theta=\Upsilon$.

The augmentation $\linlen{\Upsilon,\Theta}$ can be defined by one formula in the
language of set theory which applies to all relations $(S,\le)$ and all
$\Upsilon,\Theta\subseteq S^S$. The following question asks whether
$\linlen{\Theta}$ is the minimum among all 
such definable $\Theta$-linear augmentations
(restricted to the case $\Upsilon=\Theta$). 

\begin{question}
\label{q-1}
Is it so that for any formula $\psi(x,y,z)$, if \emph{for every relation $(S,\le)$
and every $\Theta\subseteq S^S$ there is a unique $x$ such that
$\psi(x,(S,\le),\Theta)$ holds, and $(S,x)$ is a $\Theta$-linear augmentation of
$(S,\le)$}, then for every relation $(S,\le)$ and $\Theta\subseteq S^S$,
$\psi(x,(S,\le),\Theta)$ implies $\aug{\linlen{\Theta}}{x}$\textup?
\end{question}

Sometimes $\linlen{\Upsilon,\Theta}$ yields a quasi/partial order augmentation.

\begin{lem}
\label{lem:lin-qo-aug}
Let $(O,\le)$ be a quasi order, and let $\Pi$
be a subsemigroup of $(\edo(O,\allowbreak\le),\circ)$. If $\le$ is a quasi
$\Pi$-linear ordering then $\linlen{\id,\Pi}$ is a quasi order augmentation
of $\le$.
\end{lem}
\begin{proof}
Reflexitivity holds because $\aug{\le}{\linlen{\id,\Pi}}$. 
Hence we need to verify transitivity, for which we use proposition \ref{pro:lin-simple}.
To this end we assume that $p\linlen{\id,\Pi}q$ but $p\nlinlen{\id,\Pi}r$,
and prove that $q\nlinlen{\id,\Pi}r$. Since $p\nlinlen{\id,\Pi}r$,
there exists $\pi\in\Pi$ such that 
\begin{equation}
\pi(r)\le\pi(p)\textrm{.}\label{eq:pi-r-le-pi-p}
\end{equation}
First we note that $q\nleq r$; otherwise, since $p\nleq r$, we would have $p\nleq q$, 
and since $\pi\in\edo(O,\le)$, $\pi(q)\le\pi(r)$, which would together 
with $p\linlen{\id,\Pi}q$ imply
$\pi(r)\nleq\pi(p)$ contradicting~\eqref{eq:pi-r-le-pi-p}.
It remains to find $\tau\in\Pi$ such that
$\tau(r)\le\tau(q)$ to prove $q\nlinlen{\id,\Pi}r$. If $p\le q$
then by~\eqref{eq:pi-r-le-pi-p}, $\pi(r)\le\pi(p)\le\pi(q)$ as needed;
otherwise, $p\linlen{\id,\Pi}q$, the fact that $\Pi$ is a subsemigroup
 and $\Pi$-linearity together imply
that there exists $\sigma\in\Pi$ such that $\sigma\pi(p)\le\sigma\pi(q)$;
and therefore $\sigma\pi(r)\le\sigma\pi(q)$ by~\eqref{eq:pi-r-le-pi-p},
which is as needed since $\sigma\pi\in\Pi$.
\end{proof}

\begin{lem}
\label{l-17}
Let $(P,\le)$ be a partial order, 
and let $\Pi$ be a subsemigroup of $(\edo(P,\allowbreak\le\nobreak),\circ)$. 
If $\le$ is a $\Pi$-linear ordering 
then $\linlen{\id,\Pi}$ is a partial order augmentation of $\le$.
\end{lem}
\begin{proof}
By lemma~\ref{lem:lin-qo-aug}, we know that $\linlen{\id,\Pi}$ is a quasi order
augmentation. Thus we only need verify antisymmetry. Suppose $p\linlen{\id,\Pi}
q$ and $q\linlen{\id,\Pi} p$. By $\Pi$-linearity, there exists $\pi\in\Pi$ such that
$\pi(p)$ is comparable with $\pi(q)$. Assume without loss of generality that
$\pi(p)\le \pi(q)$. Then $q\linlen{\id,\Pi} p$ implies $q\le p$. But this implies
$\tau(q)\le\tau(p)$ for all $\tau\in\Pi$, and hence for some
$\tau\in\Pi\ne\emptyset$. Therefore, $p\linlen{\id,\Pi} q$ implies $p\le q$, and now
$p=q$ since $\le$ is antisymmetric.
\end{proof}

\begin{example}
\label{x-9}
In example~\ref{x-18} we observed that $(\irr,\le)$ is $\proj$-linear. Thus
$\same{\linlen{\proj}}{\le}$ by proposition~\ref{p-23}. And moreover
by lemma~\ref{l-17}, $(\irr,\linlen{\id,\proj})$
is an $(\id,\proj)$-linear partial order augmentation of $(\irr,\le)$. Indeed for
all $x,y\in\irr$, 
\begin{equation}
\label{eq:1}
x\linlen{\id,\proj}y\Iff x\le y\OR x\llfnt y,
\end{equation}
where we let $\llfnt$ denote the relation $x\llfnt y$ if $x(n)<y(n)$ for all but
finitely many~$n$ (this relation is often written $\lfnt$ in the
literature; however, this disagrees with the notation of this paper; 
see example~\ref{x-1}).
It follows that  $x\linln{\id,\proj}y$ iff $x<y$ or $x\llfnt y$.
\end{example}

\begin{example}
\label{x-10}
By lemma~\ref{l-16} and proposition~\ref{p-1}, $(\irr,\linlen{\proj,\id})$ is also an
$(\id,\proj)$-linear augmentation of $(\irr,\le)$. For all $x,y\in\irr$,
\begin{equation}
  \label{eq:45}
  x\linlen{\proj,\id} y\Iff x\le y\OR (y\llfnt x\and y\nleq x). 
\end{equation}
However, this is not even a transitive relation. 
\end{example}

\begin{example}
\label{x-19}
We consider the $\Pi_0$-linear augmentation of $(\irr,\le)$, where $\Pi_0$ was
defined in example~\ref{x-20}. Here
\begin{equation}
\begin{alignedat}{1}
  \label{eq:55}
  x\linlen{\Pi_0} y\Iff x\le y\OR
  \bigl((\pi_\even(y)\llfnt \pi_\even(x)
  &\Or\pi_\odd(x)\nleq \pi_\odd(y))\\
  \and(\pi_\even(x)\llfnt \pi_\even(y)
  &\Or\pi_\odd(y)\nleq \pi_\odd(x))\bigr).
\end{alignedat}
\end{equation}
Note that e.g.~$\pi_\even(x)\llfnt\pi_\even(y)$ iff $x(2n)<y(2n)$ for all but
finitely many~$n$.
\end{example}

\begin{example}
\label{x-22}
Let $\Pi_1$ be the family of all members of $\proj$ of the form $\pi_h$ where either
$h=\id$ or $\ran(h)\subseteq\even$. Then $\Pi_1$ is a submonoid of
$\edo(\irr,\le)$. Clearly $(\irr,\le)$ is $\Pi_1$-linear.
And by proposition~\ref{pro:lin-simple}, we have that for all $x,y\in\irr$,
\begin{equation}
  \label{eq:57}
  \begin{split}
  x\linlen{\id,\Pi_1}y&\Iff x\le y\OR \sigma(y)\nleq\sigma(x)\text{ for all }\sigma\in\Pi_1\\
  &\Iff x\le y\OR \pi_\even(x)\llfnt\pi_\even(y). 
  \end{split}
\end{equation}
\end{example}

\subsection(Strict linearity){Strict $(\Upsilon,\Theta)$-linearity}
\label{sec:strict-upsilon-theta}

In some ways the following notion is more natural than the preceding
one. For example, while we are most interested in the case $\Pi\subseteq\epi(O,\le)$,
if $\Pi$ does contain a constant function then $(\Upsilon,\Pi)$-linearity becomes
a triviality. 

\begin{defn}
We say that $\le$ is \emph{strictly} $(\Upsilon,\Theta)$\emph{-linear}
if for all $p,q\in S$,\[
\sigma(p)\le\sigma(q)\textrm{ for some $\sigma\in\Upsilon$}\OR\tau(q)<\tau(p)
\textrm{ for some $\tau\in\Theta$.}\]
When $\Upsilon=\Theta$ we say that $\le$ is \emph{strictly} $\Theta$\emph{-linear.}
\end{defn}

Strict linearity is stronger than linearity (proposition~\ref{p-26});
however, strengthening from linearity to  strict linearity is only interesting 
when $\Upsilon\ne\Theta$, as indicated in proposition~\ref{p-28}.

\begin{prop}
\label{p-43}
$\le$ is strictly $(\id,\Theta)$-linear iff for all $p,q\in S$,\[
p\nleq q\impls\tau(q)<\tau(p)\textup{ for some $\tau\in\Theta$.}\]
\end{prop}


\begin{prop}
\label{p-42}
If $\Upsilon\subseteq\Upsilon'$ and $\Theta\subseteq\Theta'$ then
strict $(\Upsilon,\Theta)$-linearity entails strict $(\Upsilon',\Theta')$-linearity.
\end{prop}

\begin{defn}
The \emph{strict} $(\Upsilon,\Theta)$\emph{-linear augmentation}
of $\le$ is defined by
\begin{equation*}
\begin{split}
p\slinlen{\Upsilon,\Theta}q\If &p\le q\OR\\ 
(&\sigma(p)\nleq\sigma(q)\textrm{ for all $\sigma\in\Upsilon$}
\and\tau(q)\nless\tau(p)\textrm{ for all $\tau\in\Theta)$.}
\end{split}
\end{equation*}
\end{defn}

\begin{lem}
\label{l-18}
Assume that there exists $\sigma\in\Upsilon$ such that
$\Upsilon\circ\sigma\subseteq\Upsilon$ and 
$\Theta\circ\sigma\subseteq\Theta$
\textup(e.g.~if $\id\in\Upsilon$\textup), and also that $\id\in\Theta$. 
Then $\slinlen{\Upsilon,\Theta}$ is a strictly 
$(\Upsilon,\Theta)$-linear augmentation of~$\le$. 
\end{lem}
\begin{proof}
Take $p,q\in S$. Since $\sigma(p)\slinlen{\Upsilon,\Theta}\sigma(q)$ witnesses
strict $(\Upsilon,\Theta)$-linearity, we assume
$\sigma(p)\nslinlen{\Upsilon,\Theta}\sigma(q)$. Then either
$\pi\sigma(p)\le\pi\sigma(q)$ for some $\pi\in\Upsilon$, or
\begin{equation}
\label{eq:47}
\tau\sigma(q)<\tau\sigma(p)
\end{equation}
for some $\tau\in\Theta$. In the former case, 
since $\pi\sigma\in\Upsilon$, this verifies strict $(\Upsilon,\Theta)$-linearity. 
In the latter case, 
we have in particular that $\tau\sigma(q)\slinlen{\Upsilon,\Theta}\tau\sigma(p)$,
and it remains to show $\tau\sigma(p)\nslinlen{\Upsilon,\Theta}\tau\sigma(q)$. 
And this is so by~\eqref{eq:47}, because $\id(\tau\sigma(q))<\id(\tau\sigma(p))$.
\end{proof}

The strict linear augmentation is in fact an augmentation of the linear augmentation
(theorem~\ref{l-24}\eqref{item:1}).

\begin{prop}
\label{pro:If-strict-lin}
If $\le$ is strictly $(\Upsilon,\Theta)$-linear
then it is equal to its strict $(\Upsilon,\Theta)$-linear augmentation.
\end{prop}

\begin{prop}
\label{p-27}
$\Upsilon\subseteq\Upsilon'$ and $\Theta\subseteq\Theta'$ imply
$\aug{\slinlen{\Upsilon',\Theta'}}{\slinlen{\Upsilon,\Theta}}$. 
\end{prop}

The main reason we found this augmentation interesting, is that in the case
$\Upsilon=\{\id\}$ it has a particularly nice form.

\begin{prop}
\label{p-25}
Let $\Pi\subseteq\edo(S,\le)$. Then
\[
p\slinlen{\id,\Pi}q\Iff\pi(q)\nless\pi(p)\textup{ for all $\pi\in\Pi$.}
\]
\end{prop}

\begin{example}
\label{x-16}
$(\irr,\linlen{\id,\proj})$ is strictly $(\id,\proj)$-linear.
In fact, in this particular example
$\same{\linlen{\id,\proj}}{\slinlen{\id,\proj}}$.
\end{example}

\begin{example}
\label{x-17}
$(\irr,\lefnt)$ is strictly $(\id,\proj)$-linear, 
while $\aug{\slinlen{\id,\proj}}{\naug{\lefnt}{\slinlen{\id,\proj}}}$.
\end{example}

\begin{example}
\label{x-23}
Let $\Pi_1$ be the submonoid of $\edo(\irr,\le)$ from example~\ref{x-22}.
By proposition~\ref{p-25}, we have that for all $x,y\in\irr$,
\begin{equation}
  \label{eq:56}
  \begin{alignedat}{2}
  x\slinlen{\id,\Pi_1}y&\Iff &&\sigma(y)\nless \sigma(x)\text{ for all }\sigma\in\Pi_1\\
  &\Iff &&y\nless x\And \\
  &&(&\pi_\even(x)\le \pi_\even(y)
  \Or \pi_\even(x)\llfnt \pi_\even(y)).
  \end{alignedat}
\end{equation}
Note that this is nontransitive.
\end{example}

\subsection(Strictness){$(\Upsilon,\Theta)$-strictness}
\label{sec:upsil-theta-strictn}

The following property states that whenever $p<q$ there is member of $\Upsilon$
which pins the strict inequality with respect to all members of $\Theta$.

\begin{defn}
We say that $\le$ is \emph{$(\parms)$-strict} if $\Upsilon$ pins the relation
$(S,<)$ with respect to $\Theta$. In the special case $\Upsilon=\Theta$ 
we say that $\le$ is $\Theta$\emph{-strict}. 
\end{defn}

\noindent Expanding the definition yields:

\begin{prop}
\label{p-9}
$\le$ is $(\Upsilon,\Theta)$-strict iff for all
$p,q\in S$,
\begin{equation*}
p<q\impls\exists\sigma\in\Upsilon\,\forall\tau\in\Theta
\spc\tau\circ\sigma(p)<\tau\circ\sigma(q)
\textrm{.}
\end{equation*}
\end{prop}

This property is of interest to us because it distinguishes between the quasi
orderings $\le$ and $\lefnt$ of $\irr$. 

\begin{example}
\label{x-12}
The quasi order $(\irr,\le)$ is not $\proj$-strict, but $(\irr,\lefnt)$ is a
$\proj$-strict quasi order.
\end{example}

In the important case where the parameters are families of endomorphisms, 
we have an equivalent formulation.

\begin{prop}
\label{p-29}
Let  $\Lambda,\Pi\subseteq\edo(S,\le)$.
Then $\le$ is $(\Lambda,\Pi)$-strict iff for all $p,q\in S$,
\begin{equation*}
  p< q\impls \exists\sigma\in\Lambda\,\forall\tau\in\Pi
  \spc \tau\circ\sigma(q)\nleq\tau\circ\sigma(p).
\end{equation*}
\end{prop}

\begin{defn}
The $(\Upsilon,\Theta)$\emph{-strictive augmentation} of $\le$ is
defined by
\begin{multline*}
p\strlen{\Upsilon,\Theta}q\If p\le q
\OR(q<p\and\forall\sigma\in\Upsilon\,\exists\tau\in\Theta
\spc\tau\circ\sigma(q)\nless\tau\circ\sigma(p))\textrm{.}
\end{multline*}
In the special case $\Upsilon=\Theta$ we refer to the $\Theta$\emph{-strictive
augmentation} of $\le$, and write~$\strlen{\Theta}$. 
\end{defn}

\noindent In the case of endomorphisms this becomes:

\begin{prop}
\label{p-30}
Let $\Lambda,\Pi\subseteq\edo(S,\le)$. Then for all $p,q\in S$,
\begin{equation*}
  p\strlen{\Lambda,\Pi} q\Iff p\le q\OR
  (q<p\and \forall\sigma\in\Lambda\,\exists\tau\in\Pi
  \spc\tau\circ\sigma(p)\le\tau\circ\sigma(q)).
\end{equation*}
\end{prop}

\begin{lem}
\label{lem:p<str-q}
Suppose that $\Theta$ is a submonoid of $(S^{S},\circ)$
and $\Theta\circ\sigma\subseteq\Theta$ for some $\sigma\in\Upsilon$
\textup(e.g.~if $\Upsilon\subseteq\Theta$\textup). Then for all $p,q\in S$, 
\begin{equation*}
\tau(p)<\tau(q)\textup{ for all $\tau\in\Theta$}\impls p\strln{\Upsilon,\Theta}q\textrm{.}
\end{equation*}
\end{lem}
\begin{proof}
Assume $\tau(p)<\tau(q)$ for all $\tau\in\Theta$. 
Since $\id\in\Theta$, $p<q$ and in particular $p\strlen{\Upsilon,\Theta}q$.
And $q\nstrlen{\Upsilon,\Theta}p$ because this is equivalent to
\begin{equation}
q\nleq p\And(p\nless
q\Or\exists\bar\sigma\in\Upsilon\,\forall\tau\in\Theta
\spc\tau\bar\sigma(p)<\tau\bar\sigma(q))\textrm{,}\label{eq:q-nleq-p-and}
\end{equation}
and for every $\tau\in\Theta$, the hypothesis implies that $\tau\sigma(p)<\tau\sigma(q)$
because $\tau\sigma\in\nobreak\Theta$. Therefore the latter clause of the
disjunction in~\eqref{eq:q-nleq-p-and} holds, and as $q\nleq p$,
equation~\eqref{eq:q-nleq-p-and} holds, 
proving that $p\strln{\Upsilon,\Theta}q$.
\end{proof}

\begin{lem}
\label{l-19}
If $\Theta\subseteq S^{S}$ is a submonoid and $\Theta\circ\sigma\subseteq\Theta$
for some $\sigma\in\Upsilon$, then $\strlen{\Upsilon,\Theta}$ is
a $(\Upsilon,\Theta)$-strict augmentation of $\le$.
\end{lem}
\begin{proof}
Suppose $p\strln{\Upsilon,\Theta}q$, i.e.~$p\strlen{\Upsilon,\Theta}q$
and $q\nstrlen{\Upsilon,\Theta}p$. Since the latter implies that
$q\nleq p$, the former implies that $p<q$. Therefore (see~\eqref{eq:q-nleq-p-and})
there exists $\bar\sigma\in\Upsilon$ such that
\begin{equation}
\label{eq:48}
\tau\bar\sigma(p)<\tau\bar\sigma(q)\espc\textrm{for all $\tau\in\Theta$.}
\end{equation}
And thus for every $\tau\in\Theta$: $\pi\tau\bar\sigma(p)<\pi\tau\bar\sigma(q)$
for all $\pi\in\Theta$, since $\pi\tau\in\Theta$; 
hence, lemma~\ref{lem:p<str-q} yields 
$\tau\bar\sigma(p)\strln{\Upsilon,\Theta}\tau\bar\sigma(q)$.
\end{proof}

In particular:

\begin{cor}
\label{o-17}
If $\Theta\subseteq S^{S}$ is a submonoid then $\strlen{\Theta}$
is a $\Theta$-strict augmentation of~$\le$.
\end{cor}

In the case of endomorphisms, the strictive augmentation is canonical,
by which we mean that it is the minimum augmentation of the given relation
that is strict for the given parameters.

\begin{lem}
\label{l-22}
Suppose that $\altaltle$ is a $(\Lambda,\Pi)$-strict augmentation of $\le$. 
Then $\aug{\strlen{\Lambda,\Pi}}{\altaltle}$. 
\end{lem}
\begin{proof}
We assume $\aug{\le}{\altaltle}$ and $\altaltle$ is $(\Lambda,\Pi)$-strict. 
Suppose $p\strlen{\Lambda,\Pi} q$. If $p\le q$ then $p\altaltle q$ as desired. 
Otherwise, by proposition~\ref{p-30}, $q<p$ and
\begin{equation}
  \label{eq:49}
  \forall \sigma\in\Lambda\,\exists\tau\in\Pi\spc\tau\sigma(p)\le\tau\sigma(q).
\end{equation}
If, by way of contradiction, $p\naltaltle q$, then $q\altaltless p$. Therefore, by
proposition~\ref{p-29}, there exists $\sigma\in\Lambda$ such that
$\tau\sigma(p)\naltaltle\tau\sigma(q)$ for all $\tau\in\Pi$. However, this is
clearly in contradiction with~\eqref{eq:49}.
\end{proof}

In any case, we at least have that:

\begin{prop}
If $\le$ is $(\Upsilon,\Theta)$-strict then it is equal to its $(\Upsilon,\Theta)$-strictive
augmentation. In particular, if $\le$ is $\Theta$-strict then it
is equal to its $\Theta$-strictive augmentation. 
\end{prop}

Since we are primarily interested in endomorphic parameters, we would also like to
have them remain endomorphisms with respect to the augmented relation. This is
indeed the case when the parameters satisfy an additional group theoretic property.

\begin{lem}
\label{p-35}
Suppose $\edoparms\subseteq\edo(S,\le)$, where $\Lambda$ is a subsemigroup and
$\Lambda\circ\Pi\subseteq\Lambda$ \textup(e.g.~if $\Lambda$ is a right ideal\textup).
Then $\edoparms\subseteq\edo(S,\strlen{\edoparms})$. 
\end{lem}
\begin{proof}
Suppose $p\strlen{\edoparms} q$. 
Take any endomorphism $\pi\in\edo(S,\le)$. 
We want to show that $\pi(p)\strlen{\edoparms}\pi(q)$. 
If $p\le q$ then $\pi(p)\le\pi(q)$, and hence $\pi(p)\strlen{\edoparms}\pi(q)$. 
Thus we can assume that $q<p$ and
\begin{equation}
\label{eq:51}
\forall\sigma\in\Lambda\,\exists\tau\in\Pi\spc\tau\sigma(p)\le\tau\sigma(q).
\end{equation} 
Similarly, we may also assume that $\pi(p)\nleq\pi(q)$,
and thus $\pi(q)<\pi(p)$ as $q<p$ implies $\pi(q)\le\pi(p)$. It now remains to show that
$\forall\sigma\in\Lambda\,\exists\tau\in\Pi\spc\tau\sigma\pi(p)\le\tau\sigma\pi(q)$.
If $\pi\in\Lambda$ then this is so by~\eqref{eq:51} as $\Lambda$ is a subsemigroup.
And if $\pi\in\Pi$ then this follows by~\eqref{eq:51} since $\sigma\pi\in\Lambda$
for all $\sigma\in\Lambda$. 
\end{proof}

\begin{cor}
\label{o-11}
If $\Pi\subseteq\edo(S,\le)$ is a subsemigroup, 
then $\Pi\subseteq\edo(S,\strlen{\Pi})$.
\end{cor}

While the strictive augmentation is canonical, it however fails in general to be a
quasi order augmentation.

\begin{example}
\label{x-15}
For all $x,y\in\irr$,
\begin{equation}
  \label{eq:4}
  \begin{alignedat}{2}
  x \strlen{\proj} y&\Iff &&x\le y\OR\\
  &&\bigl(&y<x\and \forall
  \pi\in\proj\exists\sigma\in\proj\spc \sigma\pi(x)\le \sigma\pi(y)\bigr)\\
  &\Iff &&x\le y\OR (y<x\and x\lefnt y)\\
  &\Iff &&x\le y\OR (y<x\and x\eqfnt y).
  \end{alignedat}
\end{equation}
This is nontransitive. E.g.~$\chi_0\strlen{\proj} \zero$ and
$\zero\strlen{\proj}\chi_1$, but $\chi_0\nstrlen{\proj}\chi_1$. 
\end{example}

\begin{example}
\label{x-21}
Consider $\Pi_0$ from example~\ref{x-20}. For all $x,y\in\irr$,
\begin{equation}
  \begin{alignedat}{1}
  x\strlen{\Pi_0}y\Iff &x\le y\OR\\
  (&y<x\and \pi_\odd(x)= \pi_\odd(y)
  \and \pi_\even(x)\eqfnt \pi_\even(y)).
  \end{alignedat}
\end{equation}
\end{example}

The remainder of this subsubsection is devoted to obtaining a strict quasi order
augmentation. 

\begin{prop}
\label{p-34}
If $p\strln{\parms}q$ then $p<q$ and there exists $\sigma\in\Upsilon$
such that $\tau\circ\sigma(p)<\tau\circ\sigma(q)$ for all $\tau\in\Theta$.
\end{prop}
\begin{proof}
$p\strln{\parms}q$ implies $p\le q$ because $q\nleq p$, and thus $p<q$. 
Now use equation~\eqref{eq:q-nleq-p-and}. 
\end{proof}

\begin{lem}
\label{l-20}
Let $(O,\le)$ be a quasi order. Suppose $\edoparms\subseteq\edo(O,\le)$ satisfy
$\Pi\circ\Lambda\subseteq\Lambda$. Then every cycle in $(O,\strlen{\edoparms})$ is
bidirectional. 
\end{lem}
\begin{proof}
Suppose $p\strln{\edoparms}q$ for some $p,q\in O$. 
We need to prove there is no chain from $q$ to $p$. 
The proof is by induction on the length $n=1,2,\dots$ of a chain 
$q=p_0\strlen{\edoparms}p_1\strlen{\edoparms}\cdots\strlen{\edoparms} p_{n-1}$, 
with the hypothesis that there exists $\sigma_n\in\Lambda$ such that
\begin{equation}
  \label{eq:50}
  \tau\sigma_n(p)<\tau\sigma_n(p_{n-1})\espc\text{for all $\tau\in\Lambda$}.
\end{equation}
For $n=1$ this holds by proposition~\ref{p-34}. 
To complete the induction, suppose $p_{n-1}\strlen{\edoparms}p_n$. 
If $p_{n-1}\le p_n$ then~\eqref{eq:50} holds because
$\tau\sigma_n(p_{n-1})\le\tau\sigma_n(p_n)$ implies
$\tau\sigma_n(p)<\tau\sigma_n(p_n)$ because $\le$ is transitive.
Otherwise, there exists $\tau_n\in\Pi$ such that
$\tau_n\sigma_n(p_{n-1})\le\tau_n\sigma_n(p_n)$. 
Since $\Pi\circ\Lambda\subseteq\Lambda$, $\sigma_{n+1}=\tau_n\sigma_n\in\Lambda$,
and it satisfies~\eqref{eq:50}. 

It remains to show there is no chain from $q$ to $p$. But given a chain from $q$ of
length~$n$, by the induction result there is a $\sigma_n\in\Lambda$
satisfying~\eqref{eq:50}. And this implies that $p_{n-1}\nstrlen{\edoparms} p$, because
it implies $p_{n-1}\nleq p$ as we are dealing with endomorphisms, 
and there is no $\tau\in\Lambda$ with
$\tau\sigma_n(p)\nless \tau\sigma_n(p_{n-1})$. 
\end{proof}

\begin{notation}
We write $\strtrnlen{\Upsilon,\Theta}$ for
$\bigl(\strlen{\Upsilon,\Theta}\bigr){}^{\trn}$ the transitive augmentation of
the $(\Upsilon,\Theta)$-strictive augmentation.
\end{notation}

\begin{lem}
\label{l-21}
Suppose $\edoparms\subseteq\edo(O,\le)$, 
where $\Lambda$ is a subsemigroup and $\Pi$ is a submonoid, 
satisfy $\Pi\circ\Lambda\circ\Pi=\Lambda$ 
and $\Pi\circ\sigma\subseteq\Pi$ for some $\sigma\in\Lambda$.
Then $\strtrnlen{\edoparms}$ is $(\edoparms)$-strict, and thus is a
$(\edoparms)$-strict quasi order augmentation of $\le$.
\end{lem}
\begin{proof}
Suppose that $p\strtrnln{\edoparms} q$. 
Then there exist $p=p_0\strlen{\edoparms} p_1\strlen{\edoparms}\cdots
\strlen{\edoparms} p_{n-1}=q$ such that $p_{i-1}\strln{\edoparms} p_{i}$
for some $i=1,\dots,n-1$. 
Therefore, lemma~\ref{l-19} applies since $\Pi\circ\sigma\subseteq\Pi$, and thus
there exists $\bar\sigma\in\Lambda$ such that 
\begin{equation}
\label{eq:13}
\tau\bar\sigma(p_{i-1})\strln{\Lambda,\Pi}\tau\bar\sigma(p_i)
\espc\text{for all $\tau\in\Pi$}. 
\end{equation}
It will suffice to show that
$\tau\bar\sigma(p)\strtrnln{\edoparms}\tau\bar\sigma(q)$ for all $\tau\in\Pi$. 
Indeed for every $\tau\in\Pi$, lemma~\ref{p-35} applies 
since $\Lambda\circ\Pi\subseteq\Lambda$, and thus
$\tau\bar\sigma(p_0)\strlen{\edoparms}
\cdots\strlen{\edoparms}\tau\bar\sigma(p_{n-1})$;
therefore, since $\Pi\circ\Lambda\subseteq\Lambda$, 
by lemma~\ref{l-20}, proposition~\ref{p-33} and~\eqref{eq:13}, 
$\tau\sigma(p_0)\strtrnln{\edoparms}\tau\sigma(p_{n-1})$ as needed. 
\end{proof}

And this moreover yields the minimum strict quasi augmentation.

\begin{cor}
\label{o-13}
Suppose $\edoparms\subseteq\edo(O,\le)$, 
where $\Lambda$ is a subsemigroup and $\Pi$ is a submonoid, 
satisfy $\Pi\circ\Lambda\circ\Pi=\Lambda$ 
and $\Pi\circ\sigma\subseteq\Pi$ for some $\sigma\in\Lambda$.
If $\altaltle$ is a quasi order augmentation of $\le$ that is $(\edoparms)$-strict,
then $\aug{\strtrnlen{\edoparms}}{\altaltle}$.
\end{cor}
\begin{proof}
By lemma~\ref{l-21}, $\strlen{\edoparms}$ is $(\edoparms)$-strict, and therefore by
lemma~\ref{l-22}, $\aug{\strlen{\edoparms}}{\altaltle}$. And then since $\altaltle$
is transitive, we obtain the desired conclusion from proposition~\ref{p-32}. 
\end{proof}

\begin{cor}
\label{o-12}
Suppose $\Pi\subseteq\edo(O,\le)$ is a submonoid. 
Then $\strtrnlen{\Pi}$ is the minimum
$\Pi$-strict quasi order augmentation of $\le$.
\end{cor}

Projective strictness combined with transitivity does indeed 
characterize eventual dominance, as we shall see in theorem~\ref{x-11}.

\subsection(Correctness){$(\Upsilon,\Theta)$-correctness}
\label{sec:upsil-theta-corr}

This property concerns pinning the negated inequality $\nleq$.

\begin{defn}
\label{d-3}
We say that $\le$ is \emph{$(\parms)$-correct} if $\Upsilon$ pins the relation
$(S,\nleq)$ with respect to $\Theta$. 
When $\Upsilon=\Theta$ we say that $\le$ is $\Theta$\emph{-correct}. 
\end{defn}

\begin{prop}
\label{p-10}
$\le$ is $(\Upsilon,\Theta)$-correct iff for all
$p,q\in S$, \[
p\nleq q\impls\exists\sigma\in\Upsilon\,\forall\tau\in\Theta
\spc\tau\circ\sigma(p)\nleq\tau\circ\sigma(q)\textrm{.}\]
\end{prop}

\begin{defn}
\label{d-4}
The $(\Upsilon,\Theta)$\emph{-corrective augmentation} of $\le$
is defined by\[
p\corlen{\Upsilon,\Theta}q\If p\le q\OR\forall\sigma\in\Upsilon\,\exists\tau\in\Theta
\spc\tau\circ\sigma(p)\le\tau\circ\sigma(q)\textrm{.}\]
In the case $\Upsilon=\Theta$ we write $\corlen{\Theta,\Theta}$ as $\corlen{\Theta}$
for the \emph{$\Theta$-corrective augmentation}. 
\end{defn}

\begin{lem}
\label{p-14}
Suppose $\Theta$ is a subsemigroup of $(S^{S},\circ)$, and there is a
$\sigma\in\Upsilon$ such that either $\Theta\circ\sigma\subseteq\Theta$ or
$\sigma\circ\Theta\subseteq\Theta$ \textup(e.g.~if $\id\in\Upsilon$
or $\Upsilon\cap\Theta\ne\emptyset$\textup), then $\corlen{\Upsilon,\Theta}$
is an $(\Upsilon,\Theta)$-correct augmentation of $\le$. 
\end{lem}
\begin{proof}
Suppose $p\ncorlen{\Upsilon,\Theta}q$.
This entails the existence of $\bar\sigma\in\Upsilon$ such that 
\begin{equation}
  \label{eq:3}
 \tau\bar\sigma(p)\nleq\tau\bar\sigma(q)\espc\text{for all $\tau\in\Theta$.}  
\end{equation}
Thus for every $\tau\in\Theta$, we moreover have 
$\bar\tau\sigma\tau\bar\sigma(p)\nleq\bar\tau\sigma\tau\bar\sigma(q)$
for all $\bar\tau\in\Theta$, because the hypothesis on $\sigma$ and the fact that
$\Theta$ is a semigroup imply that $\bar\tau\sigma\tau\in\Theta$.
Therefore $\tau\bar\sigma(p)\ncorlen{\Upsilon,\Theta}\tau\bar\sigma(q)$ for all $\tau\in\Theta$. 
\end{proof}

In particular:

\begin{cor}
\label{p-38}
$\corlen{\Theta}$ is a $\Theta$-correct augmentation of $\le$, 
whenever $\Theta$ is a subsemigroup.
\end{cor}

The following shows that $\corlen{\Upsilon,\Theta}$ is the minimum augmentation with
the required correctness, and thus is indeed the canonical corrective augmentation.
Notice that no assumptions are needed here on $\Upsilon$ and $\Theta$.

\begin{lem}
\label{l-30}
If $\altaltle$ is an augmentation of $\le$ that is $(\parms)$-correct, then
$\aug{\corlen{\parms}}{\altaltle}$. 
\end{lem}
\begin{proof}
Let $\altaltle$ be an augmentation of $\le$. 
Suppose that $\altaltle$ is not an augmentation of $\corlen{\Upsilon,\Theta}$, 
say $p\corlen{\Upsilon,\Theta} q$ but $p\naltaltle q$.
Then since $\aug{\le}{\altaltle}$, we have
\begin{equation}
  \label{eq:21}
  \forall\sigma\in\Upsilon\,\exists\tau\in\Theta\spc\tau\circ\sigma(p)\le\tau\circ\sigma(q),
\end{equation}
which implies that
$\forall\sigma\in\Upsilon\,\exists\tau\in\Theta
\spc\tau\circ\sigma(p)\altaltle\tau\circ\sigma(q)$. 
Thus $\altaltle$ is not $(\Upsilon,\Theta)$-correct.
\end{proof}

\begin{cor}
\label{l-11}
If $\corlen{\Upsilon,\Theta}$ is $(\Upsilon,\Theta)$-correct, then it is the minimum
augmentation of $\le$ that is $(\Upsilon,\Theta)$-correct.
\end{cor}

\begin{cor}
\label{o-14}
If $\le$ is $(\Upsilon,\Theta)$-correct then it is equal to its
$(\Upsilon,\Theta)$-correct augmentation.
\end{cor}

Normally it has a simpler form.

\begin{prop}
\label{pro:endo-cor-refine}
Let $\edoparms\subseteq\edo(S,\le)$.
Then for all $p,q\in S$, 
\[
p\corlen{\edoparms}q\Iff\forall\sigma\in\Lambda\,\exists\tau\in\Pi
\spc\tau\circ\sigma(p)\le\tau\circ\sigma(q)\textrm{.}\]
\end{prop}

Observe that all members of a subsemigroup $\Pi$ remain endomorphisms
of the augmentation $\corlen{\Pi}$, as is desired of an `endomorphic augmentation'. 
More generally:

\begin{prop}
\label{pro:cor-preserve}
Let $\edoparms\subseteq\edo(S,\le)$. Then\textup{:
 \begin{enumerate}[(a)]
 \item \textit{If $\Lambda$ is a subsemigroup then
   $\Lambda\subseteq\edo(S,\corlen{\edoparms})$.}
 \item \textit{If $\Lambda\circ\Pi\subseteq\Lambda$ then
     $\Pi\subseteq\edo(S,\corlen{\edoparms})$.}
 \end{enumerate}
}
\end{prop}
\begin{proof}
Use proposition~\ref{pro:endo-cor-refine}. 
\end{proof}

The corrective augmentation is especially nice because it is transitive.

\begin{lem}
\label{lem:cor-quasi-order}
Let $(O,\le)$ be a quasi order, and $\edoparms\subseteq\edo(O,\le)$. 
Suppose $\Lambda$ satisfies $\Pi\circ\Lambda\subseteq\Lambda$, 
and $\Pi$ is a subsemigroup. 
Then $\corlen{\Lambda,\Pi}$ is a quasi order augmentation of $\le$. 
\end{lem}
\begin{proof}
We need to verify transitivity, 
and we do so using proposition~\ref{pro:endo-cor-refine}.
Suppose $p\corlen{\edoparms}q$ and $q\corlen{\edoparms}r$. Take 
$\sigma\in\Lambda$. Then there exists $\tau\in\Pi$ 
such that $\tau\sigma(p)\le\tau\sigma(q)$,
and since $\tau\sigma\in\Lambda$ there exists $\pi\in\Pi$ such that
$\pi\tau\sigma(q)\le\pi\tau\sigma(r)$. Now we have $\pi\tau\in\Pi$,
and since $\pi\in\edo(O,\le)$ and $\le$ is transitive, 
$(\pi\tau)\sigma(p)\le(\pi\tau)\sigma(q)\le(\pi\tau)\sigma(r)$, proving that
$p\corlen{\edoparms}r$. 
\end{proof}

\begin{cor}
\label{o-15}
If $\Pi\subseteq\edo(O,\le)$ is a subsemigroup of $(\edo(O),\circ)$, then the
$\Pi$-corrective augmentation $\corlen{\Pi}$ 
is the minimum quasi order augmentation of $\le$ that is $\Pi$-correct.
\end{cor}
\begin{proof}
By corollary~\ref{p-38}, corollary~\ref{l-11} and lemma~\ref{lem:cor-quasi-order}. 
\end{proof}

\begin{example}
\label{x-6}
$(\irr,\lefnt)$ is $\proj$-correct. This is because for any $x,y\in\irr$, 
if $x\nlefnt y$ then letting $a\in\Fin^+$ be the (infinite) set of coordinates where
$x(n)>y(n)$, $\pi_{a}(x)(n)>\pi_{a}(y)(n)$ for all $n\in\N$, and thus for any
$h:\N\to\N$, $\pi_h\pi_a(x)>\pi_h\pi_a(y)$ and in particular,  $\pi_h\pi_a(x)\nleq
\pi_h\pi_a(y)$. 
\end{example}

What is more, $\proj$-correctness characterizes $\lefnt$ in terms of $\le$.

\begin{thm}
\label{u-6}
$(\irr,\lefnt)$ is the $\proj$-corrective augmentation of $(\irr,\le)$\textup;
symbolically, $\same{\lefnt}{\corlen{\proj}}$.
\end{thm}
\begin{proof}
By lemma~\ref{l-30} and example~\ref{x-6}, 
$\corlen\proj$ is a diminishment of $\lefnt$. Conversely, 
suppose that $x\lefnt y$. Then for any injection $f:\N\to\N$, as long as the range of
$g:\N\to\N$ avoids the appropriate finite set,
$\pi_g\circ\pi_f(x)\le\pi_g\circ\pi_f(y)$, proving $x\corlen{\proj}y$. 
\end{proof}

If $\Theta$ is too large then the corrective augmentation is the
complete quasi order (see section~\ref{sec:terminology}).

\begin{prop}
\label{pro:complete}
Let $\Theta\cap\edo(S,\le)\ne\emptyset$. Suppose
that every $p\nleq q$ in $S$ has a $\tau\in\Theta$ such that $\tau(p)\le\tau(q)$.
Then $\corlen{\Upsilon,\Theta}$ is the complete quasi ordering.
\end{prop}
\begin{proof}
Take $\pi\in\Theta\cap\edo(S,\le)$, and arbitrary $p,q\in S$. For
all $\sigma\in\Upsilon$, if $\sigma(p)\le\sigma(q)$ then $\pi\sigma(p)\le\pi\sigma(q)$,
and otherwise the hypothesis implies $\tau\sigma(p)\le\tau\sigma(q)$ for some $\tau\in\Theta$,
proving $p\corlen{\Upsilon,\Theta}q$.
\end{proof}

\begin{example}
\label{x-2}
The corrective augmentation $\corlen{\edo(\irr)}$ of $(\irr,\le)$ by its family of
endomorphisms is the complete quasi order. 
\end{example}

On the other hand, if $\Theta$ is too restrictive then the augmentation is trivial.

\begin{lem}
\label{l-1}
Let $(P,\le)$ be a poset. 
Suppose that every $p\nleq q$ in $P$ has a $\sigma\in\Upsilon$ such that
$\sigma(p)>\sigma(q)$, and that $\Pi\subseteq\mono(P,\le)$. Then
$\corlen{\Upsilon,\Pi}$ is the trivial augmentation,
i.e.~$\same{\corlen{\Upsilon,\Pi}}{\le}$. 
\end{lem}
\begin{proof}
Supposing $p\nleq q$, let $\sigma\in\Upsilon$ satisfy $\sigma(p)>\sigma(q)$.
Then for all $\tau\in\Pi$, $\tau\sigma(p)>\tau\sigma(q)$
because $\tau\in\mono(P,\le)$ and $\le$ is antisymmetric. Thus
$p\ncorlen{\Upsilon,\Pi}q$ by proposition~\ref{pro:endo-cor-refine}. 
\end{proof}

\begin{lem}
\label{p-15}
If $\Lambda$ and $\Pi$ both consist entirely of \emph{order reflecting}
\textup(i.e.~$\sigma(p)\le\sigma(q)$ implies $p\le q$\textup) members of
$\edo(S,\le)$, 
then $\corlen{\edoparms}$ is the trivial augmentation. 
\end{lem}
\begin{proof}
First note that the order reflecting endomorphisms form a subsemigroup 
of $(\edo(S,\le),\circ)$. 
Then note that if
$p\nleq q$ and $\sigma$ is order reflecting, then $\sigma(p)\nleq\sigma(q)$ 
by order reflection. Now proposition~\ref{pro:endo-cor-refine} applies.
\end{proof}

\subsection(Negative strictness){Negative $(\Upsilon,\Theta)$-strictness}
\label{sec:negat-upsil-theta}

Now we pin the negated inequality $\nless$. 

\begin{defn}
We say that $\le$ is \emph{negatively $(\parms)$-strict} if $\Upsilon$ pins the
relation $(S,\nless)$ with respect to $\Theta$. 
When $\Upsilon=\Theta$ we say that $\le$ is \emph{negatively $\Theta$-strict}.
\end{defn}

\begin{prop}
\label{p-11}
$\le$ is negatively $(\Upsilon,\Theta)$-strict iff for all $p,q\in S$,
\begin{equation*}
  p\nless q\impls \exists\sigma\in\Upsilon\,\forall\tau\in\Theta
  \spc\tau\circ\sigma(p)\nless\tau\circ\sigma(q).
\end{equation*}
\end{prop}

\begin{prop}
\label{p-44}
Let $\edoparms\subseteq\edo(S,\le)$. Then negative $(\edoparms)$-strictness is
equivalent to\textup: for all $p,q\in S$,
\begin{equation*}
  p\nleq q\impls \exists\sigma\in\Upsilon\,\forall\tau\in\Theta
  \spc\tau\circ\sigma(p)\nless\tau\circ\sigma(q).
\end{equation*}
\end{prop}

\begin{example}
\label{x-13}
$(\irr,\lefnt)$ is negatively $\proj$-strict, but $(\irr,\le)$ is not. 
\end{example}

\begin{defn}
The \emph{negative $(\parms)$-strictive augmentation} is defined by
\begin{equation*}
  p\negstrlen{\parms} q\If p\le q
  \OR \forall\sigma\in\Upsilon\,\exists\tau\in\Theta
  \spc\tau\circ\sigma(p)<\tau\circ\sigma(q).
\end{equation*}
Accordingly, we write $\negstrlen{\Theta}$ for $\negstrlen{\parms}$, 
the \emph{negative $\Theta$-strictive augmentation}. 
\end{defn}

\begin{lem}
\label{l-26}
Suppose $\edoparms\subseteq\edo(S,\le)$, 
where $\Lambda$ is a subsemigroup and $\Lambda\circ\Pi\subseteq\Lambda$. 
Then $\Lambda,\Pi\subseteq\edo(S,\negstrlen{\edoparms})$. 
\end{lem}
\begin{proof}
Suppose $p\negstrlen{\edoparms}q$. Take $\pi\in\edo(S,\le)$. If $p\le q$ then
clearly $\pi(p)\negstrlen{\edoparms}\pi(q)$; hence, we can assume that
$\forall\sigma\in\Lambda\,\exists\tau\in\Pi\spc\tau\sigma(p)<\tau\sigma(q)$. 
If $\pi\in\Lambda\cup\Pi$ then
$\forall\sigma\in\Lambda\,\exists\tau\in\Pi\spc\tau\sigma\pi(p)<\tau\sigma\pi(q)$
since $\Lambda\circ(\Lambda\cup\Pi)\subseteq\Lambda$, proving that
$\pi(p)\negstrlen{\edoparms}\pi(q)$.
\end{proof}

\begin{cor}
\label{o-16}
If $\Pi\subseteq\edo(S,\le)$ is a subsemigroup then
$\Pi\subseteq\edo(S,\negstrlen{\Pi})$. 
\end{cor}

\begin{prop}
\label{l-25}
Suppose $\Theta\circ\sigma\subseteq\Theta$ for some $\sigma\in\Upsilon$. 
Then for all $p,q\in S$,
  \begin{equation*}
    \tau(p)\nless\tau(q)\textup{ for all $\tau\in\Theta$}
    \impls p\nnegstrlen{\parms} q\Or p\le q.
  \end{equation*}
\end{prop}
\begin{proof}
Suppose that $p\negstrlen{\parms}q$ and $p\nleq q$. 
This implies that there exists $\tau\in\Theta$
such that $\tau\sigma(p)<\tau\sigma(q)$. 
And since $\tau\sigma\in\Theta$, the proof is complete.
\end{proof}

\begin{lem}
\label{l-27}
Suppose $\edoparms\subseteq\edo(S,\le)$ are both subsemigroups,
$\Pi\circ\sigma\subseteq\Pi$ for some $\sigma\in\Lambda$, 
and $\Lambda\circ\Pi\subseteq\Lambda$. 
Then $\negstrlen{\edoparms}$ is negatively $(\edoparms)$-strict.
\end{lem}
\begin{proof}
Note that by lemma~\ref{l-26}, proposition~\ref{p-44} applies. Thus it suffices to
assume that $p\nnegstrlen{\edoparms} q$, and prove that there exists
$\bar\sigma\in\Lambda$ such that
$\tau\bar\sigma(p)\nnegstrln{\edoparms}\tau\bar\sigma(q)$ for all $\tau\in\Pi$. 

By this assumption,  there exists
$\bar\sigma\in\Lambda$ such that $\tau\bar\sigma(p)\nless \tau\bar\sigma(q)$ for all
$\tau\in\Pi$. Thus for every $\tau\in\Pi$, 
$\pi\tau\bar\sigma(p)\nless\pi\tau\bar\sigma(q)$ 
for all $\pi\in\Pi$ since $\pi\tau\in\Pi$, 
and thus either $\tau\bar\sigma(p)\nnegstrlen{\edoparms}\tau\bar\sigma(q)$ or else
$\tau\bar\sigma(p)\le\tau\bar\sigma(q)$, by proposition~\ref{l-25}. In the former
case our goal has been achieved, and in the latter case we obtain
$\tau\bar\sigma(q)\le\tau\bar\sigma(p)$ implying that
$\tau\bar\sigma(q)\negstrlen{\edoparms}\tau\bar\sigma(p)$ which suffices.
\end{proof}

\begin{cor}
\label{o-18}
If $\Pi\subseteq\edo(S,\le)$ is a subsemigroup then $\negstrlen{\Pi}$ is a
negatively $\Pi$-strict augmentation of $\le$.
\end{cor}

It appears that one does not get as canonical an augmentation as for the previous 
two cases of pinning. In lemma~\ref{l-29} below 
it is shown that the negatively strictive augmentation is minimal; 
however, it does not seem to be the minimum negatively strict augmentation. 
We will not attempt to confirm this with a counterexample. 

\begin{lem}
\label{l-28}
Suppose $\edoparms\subseteq\edo(S,\le)$ 
and $\Pi\circ\sigma\subseteq\Lambda$ for
some $\sigma\in\Lambda$. Then for all $p,q\in S$,
\begin{equation*}
  p\negstreq{\edoparms}q\Iff[p]=[q],
\end{equation*}
i.e.~$[p]=[q]$ means $p\le q$ and $q\le p$. 
\end{lem}
\begin{proof}
Suppose $p\negstreq{\parms}q$, but that by way of contradiction, $p\nleq q$. 
Then 
\begin{equation}
  \label{eq:53}
  \tau\sigma(p)<\tau\sigma(q)
\end{equation}
for some $\tau\in\Pi$. However, this implies that $q\nleq p$ 
because $\tau\sigma\in\edo(S)$. But since $\tau\sigma\in\Lambda$,
there must exist $\pi\in\Pi$ such that $\pi\tau\sigma(q)<\pi\tau\sigma(p)$,
contradicting~\eqref{eq:53}. The converse holds because
$\aug{\le}{\negstrlen{\edoparms}}$. 
\end{proof}

\begin{lem}
\label{l-29}
Suppose $\edoparms\subseteq\edo(S,\le)$ 
and $\Pi\circ\sigma\subseteq\Lambda$ for
some $\sigma\in\Lambda$. If $\altaltle$ is an augmentation of $\le$ that is
negatively $(\edoparms)$-strict, and moreover $\altaltle$ is a diminishment of
$\negstrlen{\edoparms}$, then in fact $\same{\altaltle}{\negstrlen{\edoparms}}$.
\end{lem}
\begin{proof}
Let $\aug{\aug{\le}{\altaltle}}{\negstrlen{\edoparms}}$ be negatively
$(\edoparms)$-strict. 
Suppose $p\negstrlen{\edoparms}q$. If $p\le q$ then $p\altaltle q$ as desired.
Assume then that $p\nleq q$, and thus
\begin{equation}
  \label{eq:52}
  \forall\bar\sigma\in\Lambda\,\exists\tau\in\Pi
  \spc\tau\bar\sigma(p)<\tau\bar\sigma(q).
\end{equation}
If it is not the case that $p\altaltle q$ then in particular, $p\naltaltless q$ and
thus there exists $\bar\sigma\in\Lambda$ such that 
$\tau\bar\sigma(p)\naltaltless\tau\bar\sigma(q)$ for all $\tau\in\Pi$; 
however, by~\eqref{eq:52}, there exists $\tau\in\Pi$ such that
\begin{equation}
\label{eq:14}
\tau\bar\sigma(p)<\tau\bar\sigma(q)
\end{equation}
which implies $\tau\bar\sigma(p)\altaltle\tau\bar\sigma(q)$,
and thus $\tau\bar\sigma(q)\altaltle\tau\bar\sigma(p)$ as well. Since
$\aug{\altaltle}{\negstrlen{\edoparms}}$, 
$\tau\bar\sigma(p)\negstreq{\edoparms}\tau\bar\sigma(q)$ and thus
$[\tau\bar\sigma(p)]=[\tau\bar\sigma(q)]$ by lemma~\ref{l-28},
contradicting~\eqref{eq:14}.
\end{proof}

\begin{example}
\label{x-14}
Consider $(\irr,\le)$. Then for all $x,y\in\irr$, $x\negstrlen{\proj}y$ iff $x\le y$
or $x\llfnt y$. Thus in this case,
$\same{\same{\linlen{\id,\proj}}{\slinlen{\id,\proj}}}{\negstrlen{\proj}}$. 
\end{example}

\begin{example}
\label{x-5}
Let $\Pi_0\subseteq\edo(\irr,\le)$ be the submonoid of example~\ref{x-20}. 
Then
\begin{equation}
  \label{eq:16}
  \begin{alignedat}{3}
  x\negstrlen{\Pi_0}y\Iff &x\le y\OR\\
  \bigl(&\pi_\even(x)\llfnt\pi_\even(y)
  &&\and \pi_\odd(x)&&\le\pi_\odd(y)\bigr)\OR\\
  \bigl(&\pi_\even(x)\lefnt\pi_\even(y)
  &&\and \pi_\odd(x)&&<\pi_\odd(y)\bigr).
  \end{alignedat}
\end{equation}
\end{example}
\section{Interrelationships}
\label{sec:interrelationships}

\begin{figure}
\includegraphics[bb=111 334 407 658]{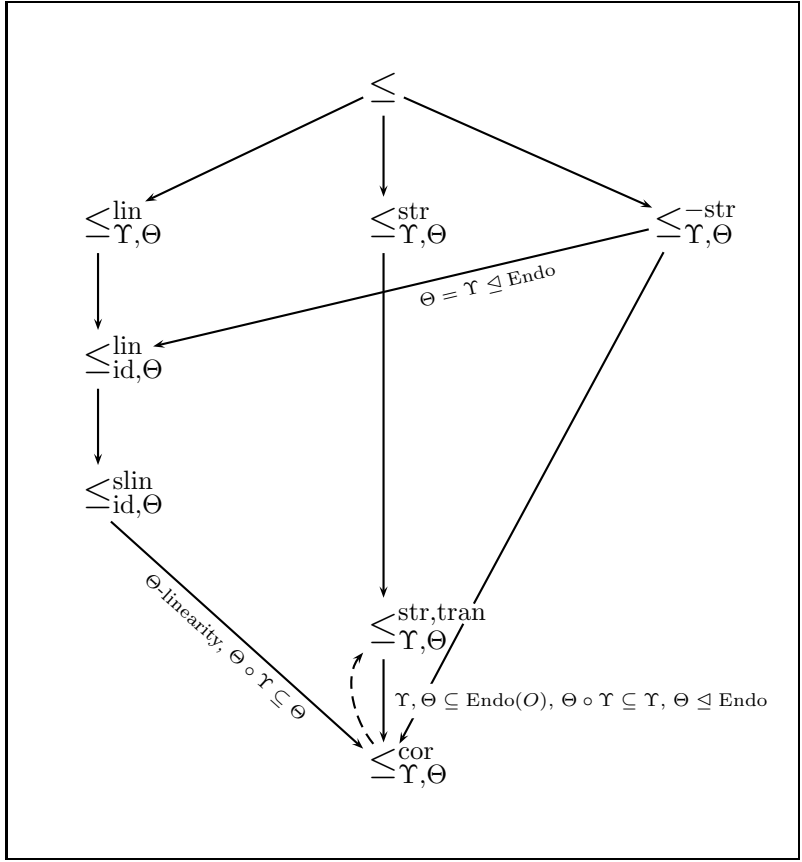}
\caption{Interrelationships}
\label{fig:1}
\end{figure}

First we examine the relationships which exist between the various properties of
relations that were introduced in section~\ref{sec:prop-quasi-orders}. Then the
corresponding relationships between the augmentations are given in
theorem~\ref{l-24}, and are summarized in figure~\ref{fig:1}. The discussion is
concluded by providing counterexamples to the other connections. 

We consider the correctness property to be especially significant, 
partly due to its position at the bottom of figure~\ref{fig:1}. 
Thus we are especially interested in theorem~\ref{u-2}, which gives a fine analysis of
the corrective augmentation as a two step augmentation (see also lemma~\ref{u-1}). For
this we need to introduce the notion of a \emph{quasi lattice}. 

\begin{defn}
A \emph{quasi lattice} is a quasi order $(L,\le)$ such that for every $p,q\in L$,
the set $\{p,q\}$ has both an infimum and supremum. We write $p\land q$ and $p\lor
q$ for the set of all infimums and supremums, respectively. The class of quasi
lattices is viewed as the category $\ql$ where the homomorphisms preserve infimums
and infimums, i.e.~$f:L\to M$ is in $\hom_{\ql}((L,\le),(M,\altle))$ iff $r\in
p\land q$ implies $f(r)\in f(p)\land f(q)$, and $s\in p\lor q$ implies $f(s)\in
f(p)\lor f(q)$, for all $p,q,r,s\in L$. 
\end{defn}

\begin{remark}
\label{r-1}
We have found various usages of the term quasi lattice in the literature.
\end{remark}

\begin{prop}
\label{p-45}
For all $p,q \in L$, $p\le q$ iff $p\in p\land q$. 
\end{prop}

\begin{prop}
\label{p-47}
$\ql$ is a subcategory of $\qo$, and thus
$\hom_{\ql}((L,\le),(M,\altle))\subseteq\hom_\relc((L,\le),(M,\altle))$. 
\end{prop}

The difference between a lattice is that there the infimum and supremum are unique
when they exist. 

\begin{prop}
\label{p-41}
$(L,\le)$ is a quasi lattice iff its antisymmetric quotient is a lattice. 
\end{prop}

\begin{example}
\label{x-3}
$(\irr,\lefnt)$ is a quasi lattice but not a lattice. And the projections
$\proj(\irr)$ are quasi lattice endomorphisms,
i.e.~$\proj(\irr)\subseteq\edo_\ql(\irr,\lefnt)$. 
\end{example}

Now we examine interrelationships between the various properties.

\begin{prop}
\label{p-26}
Strict $(\Upsilon,\Theta)$-linearity entails $(\Upsilon,\Theta)$-linearity.
\end{prop}

\begin{prop}
\label{p-28}
Strict $\Theta$-linearity is equivalent to $\Theta$-linearity.
\end{prop}

\begin{lem}
\label{l-23}
Suppose $\Theta\circ\Upsilon\subseteq\Theta$.
If $(S,\le)$ is $\Theta$-linear, then
$(\parms)$-correctness entails strict $(\id,\Theta)$-linearity.
\end{lem}
\begin{proof}
Assume that $(S,\le)$ is $(\parms)$-correct. To prove strict $(\id,\Theta)$-linearity,
suppose $p\nleq q$. By correctness, there exists $\sigma\in\Upsilon$ such that 
$\tau\sigma(p)\nleq \tau\sigma(q)$ for all $\tau\in\Theta$. And by
$\Theta$-linearity, there exists $\tau\in\Theta$ such that $\tau\sigma(p)$ is
comparable with $\tau\sigma(q)$. We conclude that $\tau\sigma(q)<\tau\sigma(p)$. 
The proof is complete with proposition~\ref{p-43}, since $\tau\sigma\in\Theta$.
\end{proof}

\begin{prop}
\label{p-40}
Suppose $\edoparms\subseteq\edo(S,\le)$. 
Then $(\edoparms)$-correctness entails $(\edoparms)$-strictness. 
\end{prop}
\begin{proof}
By proposition~\ref{p-29}. 
\end{proof}

\begin{lem}
\label{u-1}
If $(L,\le)$ is a quasi lattice and $\edoparms\subseteq\edo_{\ql}(L,\le)$, then
$(\edoparms)$-strictness implies $(\edoparms)$-correctness.
\end{lem}
\begin{proof}
Suppose $p\nleq q$ in $L$, and take $r\in p\land q$. 
Then $r<p$, and thus strictness in particular implies the
existence of $\sigma\in\Lambda$ such that 
\begin{equation}
\label{eq:54}
\tau\sigma(p)\nleq\tau\sigma(r)\espc\text{ for all $\tau\in\Pi$}.
\end{equation}
But $\tau\sigma(r)\in\tau\sigma(p)\land\tau\sigma(q)$, and thus~\eqref{eq:54}
is equivalent to $\tau\sigma(p)\nleq\tau\sigma(q)$ by proposition~\ref{p-45}, as needed.
\end{proof}

\begin{prop}
\label{p-39}
Suppose $\edoparms\subseteq\edo(S,\le)$.
Then $(\edoparms)$-correctness entails negative $(\edoparms)$-strictness.
\end{prop}
\begin{proof}
Assume $\le$ is $(\edoparms)$-correct. Suppose $p\nless q$. If moreover $p\nleq q$
then correctness gives $\sigma\in\Lambda$ 
such that $\tau\sigma(p)\nleq\tau\sigma(q)$, 
which implies $\tau\sigma(p)\nless\tau\sigma(q)$, for all $\tau\in\Pi$. 
Otherwise, $[p]=[q]$ and then since we are dealing with endomorphisms,  
for all $\sigma\in\Lambda$ and all $\tau\in\Pi$, 
$[\tau\sigma(p)]=[\tau\sigma(q)]$ which implies the desired result.
\end{proof}

\begin{prop}
\label{p-46}
Suppose $\Pi\subseteq\edo(S,\le)$ is a subsemigroup. Then $(\id,\Pi)$-linearity
entails negative $\Pi$-strictness. 
\end{prop}
\begin{proof}
By propositions~\ref{p-13} and~\ref{p-44}.
\end{proof}

\begin{thm}
\label{l-24}
We have the following relationships between various augmentations.\textup{
 \begin{enumerate}[(a)]
 \item\label{item:1} \textit{For $\Upsilon,\Theta\subseteq S^{S}$ arbitrary, 
 $\aug{\aug{\aug{\le}{\linlen{\parms}}}{\linlen{\id,\Theta}}}{\slinlen{\id,\Theta}}$.}
 \item\label{item:12} \textit{For all $\Theta\subseteq S^S$, 
 $\same{\linlen{\Theta}}{\slinlen{\Theta}}$.}
 \item\label{item:10} \textit{For all  $\parms\subseteq S^S$ where $\le$ is
     $\Theta$-linear and $\Theta\circ\Upsilon\subseteq\Theta$,
     $\aug{\slinlen{\id,\Theta}}{\corlen{\parms}}$.}
 \item\label{item:13} \textit{For all $\edoparms\subseteq\edo(S,\le)$, 
     $\aug{\strlen{\edoparms}}{\corlen{\edoparms}}$.}
 \item\label{item:11} \textit{For a quasi order $(O,\le)$, and all
     $\edoparms\subseteq\edo(O,\le)$, if $\Pi\circ\Lambda\subseteq\Lambda$
     and $\Pi$ is a subsemigroup, then $\aug{\strtrnlen{\edoparms}}{\corlen{\edoparms}}$.}
 \item\label{item:14} \textit{For $\parms\subseteq S^S$ arbitrary,
     $\aug{\negstrlen{\parms}}{\corlen{\parms}}$.}
 \item\label{item:15} \textit{For every subsemigroup $\Pi\subseteq\edo(S,\le)$, 
     $\aug{\negstrlen{\Pi}}{\linlen{\id,\Pi}}$.}
 \end{enumerate}
}
\end{thm}
\begin{proof}
\eqref{item:1} is by proposition~\ref{p-24}, 
and~\eqref{item:12} is immediate from the definitions. 

To prove~\eqref{item:10}, suppose $p\slinlen{\id,\Theta}q$. 
We need only deal with the case $p\nleq q$.
Take $\sigma\in\Upsilon$. By $\Theta$-linearity there exists
$\tau\in\Theta$ such that $\tau\sigma(p)$ is comparable to $\tau\sigma(q)$.
And even if $\tau\sigma(q)\le\tau\sigma(p)$, 
we must have $\tau\sigma(p)\le\tau\sigma(q)$ proving $p\corlen{\parms}q$,
because $\tau\sigma\in\Theta$ and thus $\tau\sigma(q)\nless\tau\sigma(p)$ 
by supposition.

\eqref{item:13} follows immediately from proposition~\ref{p-30}.
And~\eqref{item:11} is an immediate consequence of~\eqref{item:13}, 
lemma~\ref{lem:cor-quasi-order} and proposition~\ref{p-32}. 

To prove~\eqref{item:14}, suppose $p\negstrlen{\parms}q$. 
We may as well assume that $p\nleq q$. 
Then $\forall\sigma\in\Upsilon\,\exists\tau\in\Theta\spc
\tau\sigma(p)<\tau\sigma(q)$, which obviously implies that $p\corlen{\parms}q$. 

To prove \eqref{item:15} use proposition~\ref{pro:lin-simple}.
\end{proof}

\begin{thm}
\label{u-2}
Let $(L,\le)$ be a quasi lattice. Then for all
$\edoparms\subseteq\edo_{\ql}(L,\le)$, if $p\corlen{\edoparms}q$ then 
there exists $r\in L$ such that $p\strlen{\edoparms}r\strlen{\edoparms}q$. In
particular, $\aug{\corlen{\edoparms}}{\strtrnlen{\edoparms}}$.
\end{thm}
\begin{proof}
Suppose $p\corlen{\edoparms}q$ and take $r\in p\land q$. We may as well assume that
$p\nleq q$, in which case $r<p$, and
\begin{equation}
  \label{eq:15}
  \forall\sigma\in\Lambda\,\exists\tau\in\Pi\spc \tau\sigma(p)\le\tau\sigma(q).
\end{equation}
For every $\pi\in\edo_{\ql}(L,\le)$, since $\pi(r)\in\pi(p)\land\pi(q)$, if
$\pi(p)\le\pi(q)$, or equivalently $\pi(p)\in\pi(p)\land\pi(q)$, then
$\pi(p)\le\pi(r)$. Thus by~\eqref{eq:15}, proposition~\ref{p-47} and proposition~\ref{p-30},
$p\strlen{\edoparms}r$. And since $r\le q$, $r\strlen{\edoparms}q$, as required.
\end{proof}

\begin{cor}
\label{o-1}
Suppose $\edoparms\subseteq\edo_{\ql}(L,\le)$, $\Pi\circ\Lambda\subseteq\Lambda$ and
$\Pi$ is a subsemigroup. Then $\same{\corlen{\edoparms}}{\strtrnlen{\edoparms}}$. 
\end{cor}

\begin{cor}
\label{o-2}
For every subsemigroup $\Pi\subseteq\edo_{\ql}(L,\le)$,
$\same{\corlen{\Pi}}{\strtrnlen{\Pi}}$. 
\end{cor}

\begin{thm}
\label{x-11}
$(\irr,\lefnt)$ is the transitive augmentation of the $\proj$-strictive augmentation
of $(\irr,\le)$\textup; symbolically, $\same{\lefnt}{\strtrnlen{\proj}}$. 
\end{thm}
\begin{proof}
By theorem~\ref{u-6} and corollary~\ref{o-2}, $\same{\lefnt}{\strtrnlen{\proj}}$. 
\end{proof}

Note that in the figure $\Upsilon\trianglelefteq\edo$ symbolizes ``$\Upsilon$ is a
subsemigroup of $\edo$''. 

No additional arrows exist in figure~\ref{fig:1}, at least without
making additional assumptions beyond those made in theorem~\ref{l-24}:

\begin{list}{$\bullet$}
{\settowidth{\labelwidth}{$\bullet$}\setlength{\leftmargin}{\labelwidth+\labelsep}}
\item First we show that there are no other arrows originating from $\strlen{\parms}$.
Examples~\ref{x-15} and~\ref{x-16} and equation~\eqref{eq:1} demonstrate that
$\naug{\strlen{\proj}}{\slinlen{\id,\proj}}$. 
And example~\ref{x-14} implies that
$\naug{\strlen{\proj}}{\negstrlen{\proj}}$.

\item To verify that no other arrows start from $\negstrlen{\parms}$, first note that
$\naug{\negstrlen{\proj}}{\strlen{\proj}}$ by example~\ref{x-14}. It now suffices to
notice that by example~\ref{x-14}, and example~\ref{x-18} and
proposition~\ref{p-23}, $\same{\naug{\negstrlen{\proj}}{\linlen{\proj}}}{\le}$.

\item We check that no other arrows emanate from $\linlen{\parms}$. 
With examples~\ref{x-19} and~\ref{x-21}, 
we see that $\naug{\linlen{\Pi_0}}{\strlen{\Pi_0}}$. And then with
example~\ref{x-5}, we see that $\naug{\linlen{\Pi_0}}{\negstrlen{\Pi_0}}$. 

\item Clearly $\naug{\linlen{\id,\proj}}{\same{\linlen{\proj}}{\le}}$.

\item By examples~\ref{x-22} and~\ref{x-23}, $\naug{\slinlen{\id,\proj}}{\linlen{\id,\proj}}$.

\item Theorem~\ref{x-11} shows that $\naug{\strtrnlen{\proj}}{\slinlen{\id,\proj}}$, 
$\naug{\strtrnlen{\proj}}{\strlen{\proj}}$ and 
$\naug{\strtrnlen{\proj}}{\negstrlen{\proj}}$.
\end{list}

\vspace{-50pt}

\input{pinning.bbl}

\end{document}

%% file: macros.tex
\newcommand{\ran}{\operatorname{ran}}
\newcommand{\proj}{\operatorname{proj}}

\newcommand{\supp}{\operatorname{supp}}

\newcommand{\power}{\P}

\newcommand{\otp}{\operatorname{otp}}

\newcommand\bigext{^\frown}
\newcommand{\djun}{\amalg}

\renewcommand{\div}{\mathbin{/}}

\newcommand{\ds}[1]{\O_{#1}}
\newcommand{\id}{\operatornamewithlimits{id}}

\renewcommand{\hom}{\operatorname{Hom}}
\newcommand{\edo}{\operatorname{Endo}}
\newcommand{\mono}{\operatorname{Mono}}

\newcommand{\epi}{\operatorname{Epi}}

\newcommand{\aut}{\operatorname{Aut}}

\newcommand{\sym}{\operatorname{sym}}

\newcommand{\compat}{\approx}
\newcommand{\incompat}{\not\approx}
\newcommand{\rel}{\mathrel{R}}

\newcommand{\@@eq}[2]{=_{#1}^{#2}}
\newcommand{\@@neq}[2]{\ne_{#1}^{#2}}
\newcommand{\@@l}[2]{<_{#1}^{#2}}
\newcommand{\@@nl}[2]{\nless_{#1}^{#2}}
\newcommand{\@@le}[2]{\le_{#1}^{#2}}
\newcommand{\@@nle}[2]{\nleq_{#1}^{#2}}
\newcommand{\@@equiv}[2]{\equiv_{#1}^{#2}}

\newcommand{\@in}[1]{\in_#1}
\newcommand{\@eq}[1]{=_#1}
\newcommand{\@le}[1]{\le_#1}
\newcommand{\@nle}[1]{\nleq_#1}
\newcommand{\@ge}[1]{\ge_#1}

\newcommand{\u@le}[1]{\le^#1}
\newcommand{\u@l}[1]{<^#1}
\newcommand{\u@eq}[1]{=^#1}

\newcommand{\@lefnt}[1]{\@@le#1*}
\newcommand{\@nlefnt}[1]{\@@nle#1*}
\newcommand{\@eqfnt}[1]{\@@eq#1*}

\newcommand{\subseteqfnt}{\subseteq^*}

\newcommand{\lefnt}{\mathrel{\le^*}}
\newcommand{\nlefnt}{\nleq^*}

\newcommand{\eqfnt}{\mathrel{=^*}}
\newcommand{\neqfnt}{\mathrel{\ne^*}}
\newcommand{\lfnt}{<^*}
\newcommand{\nlfnt}{\nless^*}

\newcommand{\len}[1]{\le_{#1}}
\newcommand{\altle}{\lesssim}

\newcommand{\altaltle}{\sqsubseteq}
\newcommand{\naltaltle}{\not\sqsubseteq}
\newcommand{\altaltless}{\sqsubset}
\newcommand{\naltaltless}{\not\sqsubset}

\newcommand{\str}{\mathrm{str}}
\newcommand{\sep}{\mathrm{sep}}
\newcommand{\asym}{\mathrm{asym}}
\newcommand{\trn}{\mathrm{tran}}
\newcommand{\negstr}{-\str}

\newcommand{\linlen}[1]{\le_{#1}^{\mathrm{lin}}}
\newcommand{\nlinlen}[1]{\nleq_{#1}^{\mathrm{lin}}}
\newcommand{\linln}[1]{<_{#1}^{\mathrm{lin}}}

\newcommand{\slinlen}[1]{\le_{#1}^{\mathrm{slin}}}
\newcommand{\nslinlen}[1]{\nleq_{#1}^{\mathrm{slin}}}
\newcommand{\strlen}[1]{\@@le{#1}\str}
\newcommand{\nstrlen}[1]{\@@nle{#1}\str}
\newcommand{\strln}[1]{\@@l{#1}\str}
\newcommand{\nstrln}[1]{\@@nl{#1}\str}
\newcommand{\strequiv}[1]{\@@equiv{#1}\str}
\newcommand{\corlen}[1]{\le_{#1}^{\mathrm{cor}}}
\newcommand{\ncorlen}[1]{\nleq_{#1}^{\mathrm{cor}}}
\newcommand{\seple}{\le^\sep}
\newcommand{\lenn}[2]{\@@le{#1}{#2}}
\newcommand{\nlenn}[2]{\@@nle{#1}{#2}}
\newcommand{\lnn}[2]{\@@l{#1}{#2}}
\newcommand{\nlnn}[2]{\@@nl{#1}{#2}}
\newcommand{\nlen}[1]{\nleq_{#1}}

\newcommand{\eqnn}[2]{\@@eq{#1}{#2}}
\newcommand{\neqnn}[2]{\@@neq{#1}{#2}}
\newcommand{\asymle}{\u@le\asym}
\newcommand{\asyml}{\u@l\asym}
\newcommand{\llfnt}{\ll^{*}}

\newcommand{\strtrnlen}[1]{\@@le{#1}{\str,\trn}}
\newcommand{\strtrnln}[1]{\@@l{#1}{\str,\trn}}
\newcommand{\trnle}{\u@le{\trn}}
\newcommand{\trnl}{\u@l{\trn}}
\newcommand{\trneq}{\u@eq{\trn}}
\newcommand{\negstrlen}[1]{\@@le{#1}\negstr}
\newcommand{\nnegstrlen}[1]{\@@nle{#1}\negstr}
\newcommand{\nnegstrln}[1]{\@@nl{#1}\negstr}
\newcommand{\negstreq}[1]{\@@eq{#1}\negstr}
\newcommand{\altaltlen}[1]{\altaltle_{#1}}

\renewcommand{\aa}{\mathrm{aa}}
\newcommand{\eqaa}{\@eq\aa}
\newcommand{\inaa}{\@in\aa}

\newcommand{\leaa}{\@le\aa}
\newcommand{\nleaa}{\@nle\aa}
\newcommand{\geaa}{\@ge\aa}
\newcommand{\lefntaa}{\@lefnt\aa}
\newcommand{\nlefntaa}{\@nlefnt\aa}
\newcommand{\eqfntaa}{\@eqfnt\aa}

\newcommand{\naa}{\mathrm{naa}}
\newcommand{\eqnaa}{\@eq\naa}
\newcommand{\eqfntnaa}{\@eqfnt\naa}
\newcommand{\lenaa}{\@le\naa}
\newcommand{\lefntnaa}{\@lefnt\naa}

\newcommand{\eqae}{\@eq\ae}
\newcommand{\leae}{\@le\ae}

\newcommand{\Iff}{\espc\mathrm{iff}\espc}
\newcommand{\If}{\espc\textrm{if}\espc}
\newcommand{\impls}{\espc\mathrm{implies}\espc}
\renewcommand{\And}{\espc\mathrm{and}\espc}

\renewcommand{\and}{\hespc\mathrm{and}\hespc}

\newcommand{\OR}{\espc\mathrm{or}\espc}

\newcommand{\Or}{\hespc\mathrm{or}\hespc}

\renewcommand{\i}{\mathfrak i}

\newcommand{\zero}{\mathbf{0}}

\newcommand{\Fin}{\mathrm{Fin}}
\newcommand{\pN}{\power(\N)}

\newcommand{\N}{\mathbb N}

\newcommand{\even}{\mathrm{even}}
\newcommand{\odd}{\mathrm{odd}}

\newcommand{\irrationals}{\mathbb N^{\hsp\mathbb N}}
\newcommand{\irri}[1]{\N^{\hsp#1}}

\newcommand{\finp}{\Fin^{+}}

\newcommand{\pnfin}{\pN\div\Fin}

\newcommand{\relc}{\mathbf{Rel}}
\newcommand{\qo}{\mathbf{Q}}
\newcommand{\po}{\mathbf{P}}

\newcommand{\ql}{\mathbf{QL}}

\DeclareFontFamily{U}{cmsy}{}
\DeclareFontShape{U}{cmsy}{m}{n}{<12> sfixed * [10] cmsy10 
<10> <9> <8> <7> <6> <5> sfixed * [10] cmsy10}{}
\DeclareSymbolFont{customtwo}{U}{cmsy}{m}{n} 
\DeclareMathSymbol{\sctn}{\mathord}{customtwo}{"78}

\DeclareFontFamily{U}{cmmi}{}
\DeclareFontShape{U}{cmmi}{m}{n}{<20> sfixed * [11] cmmib10 <12> sfixed * [10]
cmmi10 <10> <9> <8> sfixed * [6] cmmi6 <5> <6> <7> sfixed * [5] cmmi5}{}
\DeclareSymbolFont{custom}{U}{cmmi}{m}{n}
\DeclareMathSymbol{\rharpoon}{\mathord}{custom}{"2A}
\newlength{\widt}
\newlength{\widttwo}
\newlength{\hgt}

\newcommand{\irr}{\irrationals}

\newcommand{\espc}{\quad}
\newlength{\@@hespc}
\newcommand{\hespc}{\hspace{\@@hespc}}

\newcommand{\Th}{{^{\mathrm{th}}}}

\renewcommand{\ae}{\mathrm{ae}}

\newcommand{\ulc}{\ulcorner}
\newcommand{\urc}{\urcorner}

\renewcommand{\O}{\mathcal O}
\renewcommand{\P}{\Pcal}
\newcommand{\Pcal}{\mathcal P}

\newcommand{\hsp}{\mspace{1.5mu}}

\newcommand{\spc}{\,\,\,}

\newcounter{saveenumi}



%% file: pinning.bbl
\providecommand{\bysame}{\leavevmode\hbox to3em{\hrulefill}\thinspace}
\providecommand{\MR}{\relax\ifhmode\unskip\space\fi MR }
\providecommand{\MRhref}[2]{%
  \href{http://www.ams.org/mathscinet-getitem?mr=#1}{#2}
}
\providecommand{\href}[2]{#2}
\renewcommand{\MR}[1]{}

%% file: pinning.bbl
\begin{thebibliography}{Kec95}

\bibitem[DP02]{MR1902334}
B.~A. Davey and H.~A. Priestley, \emph{Introduction to lattices and order},
  second ed., Cambridge University Press, New York, 2002. \MR{MR1902334
  (2003e:06001)}

\bibitem[Fra00]{MR1808172}
Roland Fra{\"{\i}}ss{\'e}, \emph{Theory of relations}, revised ed., Studies in
  Logic and the Foundations of Mathematics, vol. 145, North-Holland Publishing
  Co., Amsterdam, 2000, With an appendix by Norbert Sauer. \MR{MR1808172
  (2002d:03084)}

\bibitem[Hir06]{irrationals}
James Hirschorn, \emph{Characterizing the ordering of the irrationals by
  eventual dominance}, preprint, 2006.

\bibitem[Kec95]{MR1321597}
Alexander~S. Kechris, \emph{Classical descriptive set theory}, Graduate Texts
  in Mathematics, vol. 156, Springer-Verlag, New York, 1995. \MR{MR1321597
  (96e:03057)}

\end{thebibliography}
